\documentclass[oneside,english,british]{amsart}
\usepackage[T1]{fontenc}
\usepackage[latin9]{inputenc}
\usepackage{verbatim}
\usepackage{amsthm}
\usepackage{amstext}
\usepackage{amssymb}

\makeatletter
\numberwithin{equation}{section}
\numberwithin{figure}{section}
  \theoremstyle{remark}
  \newtheorem*{acknowledgement*}{Acknowledgement}
\theoremstyle{plain}
\newtheorem{thm}{Theorem}[section]
  \theoremstyle{plain}
  \newtheorem{lem}[thm]{Lemma}
  \theoremstyle{remark}
  \newtheorem*{rem*}{Remark}
  \theoremstyle{plain}
  \newtheorem{prop}[thm]{Proposition}
  \theoremstyle{plain}
  \newtheorem{cor}[thm]{Corollary}
  \theoremstyle{definition}
  \newtheorem{defn}[thm]{Definition}

\usepackage[all]{xy}
\usepackage{mathrsfs}
\usepackage{alex}

\newtheorem*{thm1}{Theorem \ref{thm:Inadequacy}}
\newtheorem*{thm2}{Theorem \ref{thm:EDL var GL2-SL2}}

\date{}

\makeatother

\usepackage{babel}

\begin{document}
\selectlanguage{english}%
\title[Extended Deligne-Lusztig varieties]{Extended Deligne-Lusztig varieties for general and special linear groups}

\author{Alexander Stasinski}
\begin{abstract}
We give a generalisation of Deligne-Lusztig varieties for general
and special linear groups over finite quotients of the ring of integers
in a non-archimedean local field. Previously, a generalisation was
given by Lusztig by attaching certain varieties to unramified maximal
tori inside Borel subgroups. In this paper we associate a family of
so-called extended Deligne-Lusztig varieties to all tamely ramified
maximal tori of the group. 

Moreover, we analyse the structure of various generalised Deligne-Lusztig
varieties, and show that the {}``unramified'' varieties, including
a certain natural generalisation, do not produce all the irreducible
representations in general. On the other hand, we prove results which
together with some computations of Lusztig show that for $\SL_{2}(\mathbb{F}_{q}[[\varpi]]/(\varpi^{2}))$,
with odd $q$, the extended Deligne-Lusztig varieties do indeed afford
all the irreducible representations.
\end{abstract}

\address{School of Mathematics\\
University of Southampton\\
Southampton, SO17 1BJ\\
UK}

\email{\texttt{a.stasinski@soton.ac.uk}}

\maketitle

\section{Introduction}

Let $F$ be a non-archimedean local field with finite residue field
$\mathbb{F}_{q}$. Let $\mathcal{O}_{F}$ be the ring of integers
in $F$, and let $\mathfrak{p}$ be its maximal ideal. If $r\geq1$
is a natural number, we write $\mathcal{O}_{F,r}$ for the finite
quotient ring $\mathcal{O}_{F}/\mathfrak{p}^{r}$. Let $\mathbf{G}$
be a reductive group scheme over $\mathcal{O}_{F}$. The representation
theory of groups of the form $\mathbf{G}(\mathcal{O}_{F,r})$, in
particular for $\mathbf{G}=\textrm{GL}_{n}$, has recently attracted
attention from several different directions. On the one hand, there
are the {}``algebraic'' approaches to the construction of representations.
These include the method of Clifford theory and conjugacy orbits,
which can deal explicitly with the class of regular representations
(cf.~\cite{Hill_regular} and \cite{Alex_smooth_reps_GL2}). Another
approach, due to Onn \cite{Uri-rank-2}, is based on a generalisation
of parabolic induction for general automorphism groups of finite $\mathcal{O}_{F}$-modules.
This approach and the associated notion of cuspidality for $\GL_{n}(\mathcal{O}_{F,r})$
are developed in \cite{AOPS}. Moreover, by the work of Henniart \cite{Henniart-appendix}
and Paskunas \cite{Vytas-unicity}, it is known that every supercuspidal
representation of $\GL_{n}(F)$ has a unique type on $\GL_{n}(\mathcal{O}_{F})$.
Hence the representation theory of the finite groups $\GL_{n}(\mathcal{O}_{F,r})$
encodes important information about the infinite-dimensional representation
theory of the $p$-adic group $\GL_{n}(F)$.

On the other hand, there is the cohomological approach to constructing
representations. The case $r=1$ corresponds to connected reductive
groups over finite fields and was treated in the celebrated work of
Deligne and Lusztig \cite{delignelusztig}. In \cite{Springer-Kloosterman},
Springer asks whether the geometric methods employed for $r=1$ can
be used to deal also with groups of the form $\mathbf{G}(\mathcal{O}_{F,r})$,
for $r\geq2$. The first step in this direction was taken by Lusztig
\cite{lusztig}, where a cohomological construction of certain representations
of groups of the form $\mathbf{G}(\mathcal{O}_{F,r})$ was suggested
(without proof). More recently, the proof was given in \cite{Lusztig-Fin-Rings}
for the case where $F$ is of positive characteristic, and this was
generalised to groups over arbitrary finite local rings in \cite{Alex_Unramified_reps}.
This construction attaches varieties and corresponding virtual representations
$R_{T,U}(\theta)$ of $\mathbf{G}(\mathcal{O}_{F,r})$ to certain
maximal tori in $\mathbf{G}$. However, this construction has two
limitations. Firstly, in contrast to the case $r=1$, it is not true
for $r\geq2$ that every irreducible representation of $\mathbf{G}(\mathcal{O}_{F,r})$
is a component of some $R_{T,U}(\theta)$. Secondly, the maximal tori
in $\mathbf{G}$ correspond to unramified tori in the group $\mathbf{G}\times F$,
that is, maximal tori which are split after an unramified extension.
However, there also exist ramified maximal tori in $\mathbf{G}\times F$,
and these are known to play a role in the representation theory of
$\GL_{n}(\mathcal{O}_{F,r})$ and $\SL_{n}(\mathcal{O}_{F,r})$ analogous
to that of the unramified maximal tori. In particular, since the work
of Howe \cite{Howe-Tame} it has been known that tamely ramified supercuspidal
representations of $\GL_{n}(F)$ come in families attached to maximal
tori. Given the correspondence between supercuspidal representations
of $\GL_{n}(F)$ and their types on $\GL_{n}(\mathcal{O}_{F})$, it
is not surprising that ramified maximal tori should play a role in
the representation theory of $\GL_{n}(\mathcal{O}_{F,r})$.

It is thus natural to ask whether it is possible to generalise the
{}``unramified'' construction of \cite{Lusztig-Fin-Rings} and \cite{Alex_Unramified_reps}
to account also for the ramified maximal tori. The main purpose of
this paper is to introduce a family of so-called \emph{extended Deligne-Lusztig
varieties}, corresponding to all the tamely ramified maximal tori.
Another part of the paper motivates our approach by showing the inadequacy
of varieties defined only with respect to unramified extensions of
$F$. Finally, we show in a non-trivial special case that our construction
leads to the expected result, namely, that varieties attached to a
ramified maximal torus realise in their cohomology a family of representations
which is known (by the algebraic construction) to be associated to
this maximal torus.

The following is a more detailed outline of the paper. For a scheme
$\mathbf{X}$ over $\overline{\mathbb{F}}_{q}$, and a prime $l$
different from $p$, we will consider the $l$-adic \'etale cohomology
groups with compact support $H_{c}^{i}(\mathbf{X},\overline{\mathbb{Q}}_{l})$.
In what follows, $l$ will be fixed and we will denote $H_{c}^{i}(\mathbf{X},\overline{\mathbb{Q}}_{l})$
simply by $H_{c}^{i}(\mathbf{X})$. We denote the alternating sum
of cohomologies $\sum_{i\geq0}(-1)^{i}H_{c}^{i}(\mathbf{X})$ by $H_{c}^{*}(\mathbf{X})$.
Let $F^{\mathrm{ur}}$ be the maximal unramified extension of $F$
(inside a fixed algebraic closure of $F$), and let $\mathcal{O}_{F^{\mathrm{ur}}}$
be its ring of integers. The construction of \cite{Lusztig-Fin-Rings}
and \cite{Alex_Unramified_reps} considers the finite group $\mathbf{G}(\mathcal{O}_{F,r})$
as the fixed-point subgroup of $\mathbf{G}(\mathcal{O}_{F^{\mathrm{ur}},r})$
under a Frobenius endomorphism $\phi:G_{r}\rightarrow G_{r}$, typically
induced by the (arithmetic) Frobenius element in $\Gal(F^{\mathrm{ur}}/F)$.
The Greenberg functor allows one to view $\mathbf{G}(\mathcal{O}_{F^{\mathrm{ur}},r})$
as a connected affine algebraic group $G_{r}$ over the algebraic
closure $\overline{\mathbb{F}}_{q}$, and $\mathbf{G}(\mathcal{O}_{F,r})$
is naturally isomorphic to a subgroup $G_{F,r}$ of $G_{r}$. For
instance, if $\phi$ comes from the Frobenius in $\Gal(F^{\mathrm{ur}}/F)$,
then $G_{r}^{\phi}\cong G_{F,r}$. Similarly, for every subgroup scheme
$\mathbf{H}$ of $\mathbf{G}$, we have a connected algebraic subgroup
$H_{r}\cong\mathbf{H}(\mathcal{O}_{F^{\mathrm{ur}},r})$ of $G_{r}$.
For $r\geq r'\geq1$ we have a natural map $\rho_{r,r'}:H_{r}\rightarrow H_{r'}$,
and we denote its kernel by~$H_{r}^{r'}$. 

Suppose that $\mathbf{T}$ is a maximal torus in $\mathbf{G}\times\mathcal{O}_{F^{\mathrm{ur}}}$
contained in a Borel subgroup $\mathbf{B}$ with unipotent radical
$\mathbf{U}$ such that $T{}_{r}$ and $U_{r}$ are $\phi$-stable.
Let $L:G_{r}\rightarrow G_{r}$ be the Lang map, given by $g\mapsto g^{-1}\phi(g)$.
For any element $w$ in the Weyl group $N_{G_{1}}(T_{1})/T_{1}$,
and any lift $\hat{w}\in N_{G_{r}}(T_{r})$ of $w$, we can then define
the varieties\begin{eqnarray*}
 &  & X_{r}(w)=L^{-1}(\dot{w}B_{r})/B_{r}\cap\dot{w}B_{r}\dot{w}^{-1},\\
 &  & \widetilde{X}_{r}(\hat{w})=L^{-1}(\hat{w}U_{r})/U_{r}\cap\hat{w}U_{r}\hat{w}^{-1},\end{eqnarray*}
where $\widetilde{X}_{r}(\hat{w})$ is a finite cover of $X_{r}(w)$.
These varieties were first considered by Lusztig \cite{lusztig},
and coincide with classical Deligne-Lusztig varieties for $r=1$.
For $r=1$ the Bruhat decomposition in $G_{1}$ implies that the varieties
$X_{1}(w)$, and hence the corresponding covers $\widetilde{X}_{1}(\hat{w})$,
are attached to double $B_{1}$-$B_{1}$ cosets. 

It was shown by Deligne and Lusztig \cite{delignelusztig} that every
irreducible representation of $G_{1}^{\phi}$ is a component of the
cohomology of some variety $\widetilde{X}_{1}(\hat{w})$. In contrast,
using the varieties $\widetilde{X}_{r}(\hat{w})$ for $r\geq2$, this
is no longer true in general. On the other hand, for $r\geq2$ there
exist double $B_{r}$-$B_{r}$ cosets which are not indexed by elements
of the Weyl group. In order to construct the missing representations
it therefore seems natural to define the following varieties (first
considered by Lusztig)\[
L^{-1}(xB_{r})/B_{r}\cap xB_{r}x^{-1},\quad L^{-1}(xU_{r})/U_{r}\cap xU_{r}x^{-1},\quad\text{for any }x\in G_{r}.\]
One may then hope that since these varieties account for all double
$B_{r}$-$B_{r}$ cosets in $G_{r}$, they may also afford further
representations of $G_{r}^{\phi}$, not obtainable by the varieties
$\widetilde{X}_{r}(\hat{w})$. However, it turns out that this is
not the case, and we prove in Section~\ref{sec:The-unramified-approach}
that there are non-trivial cases where these varieties do not afford
any new representations beyond those given by the varieties $\widetilde{X}_{r}(\hat{w})$.
In Subsection~\ref{sub:Reps-SL_2} we give an explicit algebraic
description of the irreducible representations of $\SL_{2}(\mathcal{O}_{F,r})$,
using Clifford theory and orbits. This construction is well-known
for odd $q$, but the case when $q$ is a power of $2$ requires a
modification and does not seem to have previously appeared in this
form. 

Assume for the moment that $\mathbf{G}=\SL_{2}$, and let $\mathbf{U}$
and $\mathbf{U}^{-}$ be the upper and lower uni-triangular subgroups,
respectively. If $G$ is a finite group acting on two varieties $X$
and $Y$, we write $X\sim Y$ if $H_{c}^{*}(X)\cong H_{c}^{*}(Y)$
as virtual $G$-representations. In Subsection~\ref{sub:A-counterexample},
we show

\begin{thm1}Let $y\in(U^{-})_{2}^{1}$\foreignlanguage{british}{.}
Then $L^{-1}(yU_{2})\sim\widetilde{X}_{2}(1)$, and hence \[
H_{c}^{*}(L^{-1}(yU_{2}))\cong\Ind_{U_{2}^{\phi}}^{G_{2}^{\phi}}\mathbf{1}\]
as $G_{2}^{\phi}$-representations.\end{thm1}\noindent Together
with Proposition~\ref{pro:{1,w}(U^-)} and Proposition~\ref{pro:Inadeq-X_2(e)}
this result implies that any irreducible representation of $\SL_{2}(\mathcal{O}_{F,2})$
which appears in the cohomology of a variety of the form $L^{-1}(xB_{2})/B_{2}\cap xB_{2}x^{-1}$,\foreignlanguage{british}{
$L^{-1}(xU_{2})/U_{2}\cap xU_{2}x^{-1}$, or $L^{-1}(xU_{2})$ already
appears in the cohomology of a variety $\widetilde{X}_{2}(\hat{w})$,
where $w$ is one of the two elements of $N_{G_{1}}(T_{1})/T_{1}$.
Combining this with results of Lusztig on the cohomology of $\widetilde{X}_{2}(\hat{w})$,
for $w\neq1$ and $F$ of positive characteristic (cf.~\cite{Lusztig-Fin-Rings},
3), we deduce as a corollary that there exist certain nilpotent representations
of }$\SL_{2}(\mathcal{O}_{F,2})$, for $F$ of positive characteristic,
which do not appear in the cohomology of any of the above varieties.\foreignlanguage{british}{
}%
{}

\selectlanguage{british}%
Having shown that the idea of attaching generalised Deligne-Lusztig
varieties to double $B_{r}$-$B_{r}$ cosets does not lead to a satisfactory
construction, we turn to another point of view. \foreignlanguage{english}{In
this paper we will primarily be concerned with the cases $\mathbf{G}=\GL_{n}$
or $\mathbf{G}=\SL_{n}$, and where $\phi$ is the standard Frobenius.
Assume now that we are in one of these cases.}

Rather than using the varieties $\widetilde{X}_{r}(\hat{w})$, the
unramified representations $R_{T,U}(\theta)$ of \cite{Lusztig-Fin-Rings}
and \cite{Alex_Unramified_reps} can also be constructed by using
another type of variety. A variety of this kind is attached to a Borel
subgroup containing certain maximal torus. Let now $\mathbf{T}$ be
any maximal torus of \foreignlanguage{english}{$\mathbf{G}\times\mathcal{O}_{F^{\mathrm{ur}}}$}
such that $T_{r}$ is $\phi$-stable. Let $\mathbf{B}$ be a Borel
subgroup containing $\mathbf{T}$, and let $\mathbf{U}$ be the unipotent
radical of $\mathbf{B}$.\foreignlanguage{english}{ One can then attach
a Deligne-Lusztig variety to the inclusion $T_{r}\subset B_{r}$.
In the case $r=1$, the group $T_{1}$ is a maximal torus of $G_{1}$,
but in general $T_{r}$ is not a maximal torus, but a Cartan subgroup
of $G_{r}$. A $\phi$-stable Cartan subgroup $T_{r}$ is the connected
centraliser of a regular semisimple element in $G_{r}^{\phi}$. This
shows the relation between regular semi\-simple elements in $G_{r}^{\phi}$
and the unramified Deligne-Lusztig construction. The work of Hill
\cite{Hill_regular} for $\GL_{n}$, and the results for $\SL_{2}$
(see Subsection~\ref{sub:Reps-SL_2}) clearly show that the regular
elements in $\mathbf{G}(\mathcal{O}_{F,r})$ and their centralisers
play an important role in the representation theory of $\mathbf{G}(\mathcal{O}_{F,r})$.
Among the elements in $\mathbf{G}(\mathcal{O}_{F^{\mathrm{ur}},r})$,
there are those with distinct eigenvalues in some extension of the
ring $\mathcal{O}_{F^{\mathrm{ur}},r}$. We call such elements, and
the corresponding elements in $G_{r}$, \emph{separable}. For $r=1$
they are precisely the regular semisimple elements, but in general
there are non-regular unipotent separable elements. The Cartan subgroups
$T_{r}$ are thus the reductions mod $\mathfrak{p}^{r}$ of the $\mathcal{O}_{F^{\mathrm{ur}}}$-points
of unramified maximal tori in $\mathbf{G}\times F^{\mathrm{ur}}$
defined over $\mathcal{O}_{F^{\mathrm{ur}}}$, and correspond to regular
semisimple elements. In addition, there exist subgroups of $\mathbf{G}(\mathcal{O}_{F^{\mathrm{ur}},r})$
which come from ramified tori, and these are the centralisers of regular
separable elements which are not semisimple. }

\selectlanguage{english}%
The idea in Section~\ref{sec:EDL-varieties} is that one should attach
generalised Deligne-Lusztig varieties not only to unramified maximal
tori, but to the centraliser of any regular separable element in $G_{r}^{\phi}$.
To achieve this, we consider an arbitrary regular separable element
$x\in G_{r}^{\phi}$, and its centraliser $C_{G_{r}}(x)$, called
a \emph{quasi-Cartan} subgroup. To generalise the unramified case,
we would also need an inclusion of $C_{G_{r}}(x)$ into a group of
the form $B_{r}$. However, one feature of general regular separable
elements is that they may not be triangulable in $G_{r}$, that is,
$x$ may not be conjugate in $G_{r}$ to any element in $B_{r}$.
This means that unlike the Cartan subgroups $T_{r}$, general quasi-Cartans
may not lie inside any conjugate of $B_{r}$. We are thus lead to
extend the base field $F$ to a ramified extension. More precisely,
in Section~\ref{sec:EDL-varieties} we show that given any element
$x\in G_{F,r'}$, for some $r'\geq1$, there exists a finite extension
$L/F^{\mathrm{ur}}$, an integer $r\geq r'$, a connected affine algebraic
group $G_{L,r}\cong\mathbf{G}(\mathcal{O}_{L,r})$, and a $\lambda\in G_{L,r}$,
\foreignlanguage{british}{such that $G_{F,r'}\subseteq G_{L,r}$ and
such that $\lambda^{-1}x\lambda\in B_{L,r}$. This implies that if
$x$ is regular separable, the}n\[
C_{G_{r}}(x)\subseteq\lambda B_{L,r}\lambda^{-1}.\]
Given a $\phi$-stable quasi-Cartan $C_{G_{r}}(x)$, and a group $\lambda B_{L,r}\lambda^{-1}$
containing it, and assuming that $L/F^{\mathrm{ur}}$ is tamely ramified,
we construct a variety $X_{L,r}^{\Sigma}(\lambda)$, where $\Sigma$
contains two endomorphisms of $G_{L,r}$ (including one Frobenius).
The variety $X_{L,r}^{\Sigma}(\lambda)$ is a subvariety of $G_{L,r}/B_{L,r}$,
which is a generalisation of the flag variety of Borel subgroups,
and is provided with an action of the finite groups of fixed points
$G_{L,r}^{\Sigma}$. When $L/F^{\mathrm{ur}}$ is tamely ramified,
we show that $G_{L,r}^{\Sigma}=G_{F,r'}$.

It is also important to define finite covers of $X_{L,r}^{\Sigma}(\lambda)$,
generalising $\widetilde{X}_{r}(\hat{w})$. However, in general there
does not seem to be any straightforward way to define such a cover
of the whole of $X_{L,r}^{\Sigma}(\lambda)$, but only of a certain
subvariety of $X_{L,r}^{\Sigma}(\lambda)$. The covers we construct
are denoted $\widetilde{X}_{L,r}^{\Sigma}(\lambda)$, and do indeed
reduce to the covers $\widetilde{X}_{r}(\hat{w})$ in the unramified
case. In particular, $\widetilde{X}_{L,r}^{\Sigma}(\lambda)$ also
carries an action of $G_{L,r}^{\Sigma}$, and a commuting action of
a finite group $S(\lambda)/S(\lambda)^{0}$. This generalises the
action of $G_{r}^{\phi}\times T_{r}^{\hat{w}\phi}$ on $\widetilde{X}_{r}(\hat{w})$.
We call the varieties $X_{L,r}^{\Sigma}(\lambda)$ and $\widetilde{X}_{L,r}^{\Sigma}(\lambda)$
\emph{extended Deligne-Lusztig varieties}.

\selectlanguage{british}%

\selectlanguage{english}%
In Section~\ref{sec:EDLforGL2-SL2} we study the extended Deligne-Lusztig
varieties for $\mathbf{G}=\GL_{2}$ and $\mathbf{G}=\SL_{2}$, with
$F$ of odd characteristic and $r=3$. In this case, only one (tamely)
ramified quadratic extension $L/F^{\mathrm{ur}}$ occurs, and we have
$G_{L,3}^{\Sigma}=G_{F,2}\cong\mathbf{G}(\mathcal{O}_{F,2})$. There
are four conjugacy classes of rational quasi-Cartan subgroups of $G_{2}$.
The two classes of Cartan subgroups give rise to the {}``unramified''
varieties $\widetilde{X}_{2}(1)$ and $\widetilde{X}_{2}(\dot{w})$,
respectively. The third class gives rise to an extended Deligne-Lusztig
variety $\widetilde{X}_{L,3}^{\Sigma}(\lambda)$, and\foreignlanguage{british}{
we show the following \begin{thm2}Let $\mathbf{Z}$ be the centre
of $\mathbf{G}$. Then there exists a $G_{L,3}^{\Sigma}$-equivariant
isomorphism\[
\widetilde{X}_{L,3}^{\Sigma}(\lambda)/(Z_{L,3}^{1})^{\phi}\cong G_{L,3}^{\Sigma}/(Z_{L,3}^{1})^{\Sigma}(U_{L,3}^{1})^{\Sigma}.\]
\end{thm2}\noindent Here $Z_{L,3}^{1}$ is the kernel of the natural
reduction map $Z_{L,3}\rightarrow Z_{L,1}$, and similarly for $U_{L,3}^{1}$.
Combining this result with results of} Lusztig \cite{Lusztig-Fin-Rings},
we can show in particular that every irreducible representation of
$\SL_{2}(\mathbb{F}_{q}[[\varpi]]/(\varpi^{2}))$, with odd $q$ appears
in the cohomology of some extended Deligne-Lusztig variety.

\selectlanguage{british}%
In the final section, we state some open problems and indicate several
directions in which our results could be taken further.
\selectlanguage{english}%
\begin{acknowledgement*}
The main parts of this work were carried out under EPSRC Grant EP/C527402.
The author also acknowledges support by EPSRC Grant EP/F044194/1 during
the preparation of the paper. It is a pleasure to thank \mbox{A.-M.}~Aubert,
B.~K\"ock, U.~Onn, and C.~Voll for valuable comments and encouraging
discussions.
\end{acknowledgement*}

\section{\label{sec:Preliminaries}Notation and general facts}

For any discrete valuation field $F$ we denote by $\mathcal{O}_{F}$
its ring of integers, by $\mathfrak{p}_{F}$ the maximal ideal of
$\mathcal{O}_{F}$, and by $k=k_{F}$ the residue field (which we
always assume to be perfect). If $r\geq1$ is a natural number, we
let $\mathcal{O}_{F,r}$ denote the quotient ring $\mathcal{O}_{F}/\mathfrak{p}_{F}^{r}$.
Throughout the paper $\varpi=\varpi_{F}$ will denote a fixed prime
element of $\mathcal{O}_{F}$.

Let $\mathbf{X}$ be a scheme of finite type over $\mathcal{O}_{F,r}$.
Greenberg \cite{greenberg1,greenberg2} has defined a functor $\mathcal{F}_{\mathcal{O}_{F,r}}$
from the category of schemes of finite type over $\mathcal{O}_{F,r}$
to the category of schemes over $k$, such that there exists a canonical
isomorphism \[
\mathbf{X}(\mathcal{O}_{F,r})\cong(\mathcal{F}_{\mathcal{O}_{F,r}}\mathbf{X})(k),\]
and such that $\mathcal{F}_{\mathcal{O}_{F,1}}=\mathcal{F}_{k}$ is
the identity functor. Moreover, Greenberg has shown that the functor
$\mathcal{F}_{\mathcal{O}_{F,r}}$ preserves schemes of finite type,
separated schemes, affine schemes, smooth schemes, open and closed
subschemes, and group schemes, over the corresponding bases, respectively.
If $\mathbf{X}$ is smooth over $\mathcal{O}_{F,r}$ and $\mathbf{X}\times k$
is reduced and irreducible, then $\mathcal{F}_{\mathcal{O}_{F,r}}\mathbf{X}$
is reduced and irreducible (\cite{greenberg2}, 2, Corollary~2).

Let $\mathbf{G}$ be an affine smooth group scheme over $\mathcal{O}_{F}$.
By definition it is then also of finite type over $\mathcal{O}_{F}$.
For any natural number $r\ge1$ we define\[
{G}_{F,r}:=\mathcal{F}_{\mathcal{O}_{F,r}}(\mathbf{G}\times_{\mathcal{O}_{F}}\mathcal{O}_{F,r})(k).\]
By the results of Greenberg, ${G}_{F,r}$ is then the $k$-points
of a smooth affine group scheme over $k$. It can thus be identified
with the $k$-points of an affine algebraic group defined over $k$.
Since $\mathbf{G}$ is smooth over $\mathcal{O}_{F}$, it follows
that for any natural numbers $r\geq r'\geq1$, the reduction map $\mathcal{O}_{F,r}\rightarrow\mathcal{O}_{F,r'}$
induces a surjective homomorphism $\rho_{r,r'}:{G}_{F,r}\rightarrow G_{F,r'}$.
The kernel of $\rho_{r,r'}$ is denoted by ${G}_{F,r}^{r'}$. The
multiplicative representatives map $k^{\times}\rightarrow\mathcal{O}_{F,r}^{\times}$
induces a section $i_{r}:G_{F,1}\rightarrow{G}_{F,r}$. In the case
where $F$ is of positive characteristic, there is an inclusion of
$k$-algebras $k\rightarrow\mathcal{O}_{F,r}$, and $i_{r}$ is an
injective homomorphism. When $F$ is of characteristic zero $i_{r}$
is not in general a homomorphism. However, if $\mathbf{G}$ is a split
torus, then $i_{r}$ is always a homomorphism, irrespective of the
characteristic of $F$. 

Following \cite{SGA3}, XIX 2.7, we call a group scheme $\mathbf{G}$
over a base scheme $\mathbf{S}$ \emph{reductive} if $\mathbf{G}$
is affine and smooth over $\mathbf{S}$, and if its geometric fibres
are connected and reductive as algebraic groups. If $\mathbf{G}$
is a reductive group scheme over $\mathbf{S}$, we will speak of maximal
tori and Borel subgroups of $\mathbf{G}$, which are also group schemes
over $\mathbf{S}$. For any Borel subgroup of $\mathbf{G}$ there
is also a well-defined unipotent radical. For these notions, see \cite{SGA3},
XXII 1.3, XIV 4.5, and XXVI 1.6, respectively. For more on reductive
group schemes, see \cite{Alex_Unramified_reps} and its references.

From now on and throughout the paper, let $F$ denote a local field
with finite residue field $\mathbb{F}_{q}$ of characteristic $p$.
We will use the same symbol $\mathfrak{p}_{F}$ to denote the maximal
ideal in $\mathcal{O}_{F}$, as well as the maximal ideal in any of
the quotients $\mathcal{O}_{F,r}$. Let $\mathbf{G}$ be a reductive
group scheme over $\mathcal{O}_{F}$. By definition, $\mathbf{G}$
is affine and smooth over $\mathcal{O}_{F}$. We fix an algebraic
closure of $F$ in which all algebraic extensions are taken. Denote
by $F^{\textrm{ur}}$ the maximal unramified extension of $F$ with
residue field $\overline{\mathbb{F}}_{q}$, an algebraic closure of
$\mathbb{F}_{q}$. Suppose that $L$ is a finite extension of $F^{\textrm{ur}}$.
Then $L$ also has residue field $\overline{\mathbb{F}}_{q}$. We
define \[
{G}_{L,r}:=(\mathbf{G}\times_{\mathcal{O}_{F}}\mathcal{O}_{L})_{L,r}=\mathcal{F}_{\mathcal{O}_{L,r}}(\mathbf{G}\times_{\mathcal{O}_{F}}\mathcal{O}_{L,r})(\overline{\mathbb{F}}_{q}).\]
Thus ${G}_{L,r}$ is an affine algebraic group over $\overline{\mathbb{F}}_{q}$.
Since $\mathbf{G}$ has connected fibres (by definition), ${G}_{L,r}$
is connected. For $F^{\mathrm{ur}}$ we will drop the subscript and
write ${G}_{r}$ for $G_{F^{\textrm{ur}},r}$, and ${G}_{r}^{r'}$
for the kernel $G_{F^{\textrm{ur}},r}^{r'}$.

If $G$ is a finite group, we denote by $\Irr(G)$ the set of irreducible
$\overline{\mathbb{Q}}_{l}$-represen\-tations of $G$. Since the
values of the characters in $\Irr(G)$ all lie in some finite extension
of $\mathbb{Q}$, there is a character preserving bijection between
$\Irr(G)$ and the set of irreducible complex representations of $G$.
For any finite group $G$ we denote its trivial representation by
$\mathbf{1}$.

$ $If $x$ is a real number, we will write $[x]$ for the largest
integer $\leq x$.

Many results about $l$-adic cohomology used in classical Deligne-Lusztig
theory are applicable also in the generalised situations we will consider,
and throughout we will assume familiarity with the results stated
in \cite{dignemichel}, 10. In what follows, all varieties will be
separated reduced schemes of finite type over $\overline{\mathbb{F}}_{q}$,
and we identify every variety with its set of $\overline{\mathbb{F}}_{q}$-points.
Suppose that $G$ is a finite group acting on a variety $X$. Then
each $g\in G$ induces an element of $\Aut_{\overline{\mathbb{Q}}_{l}}(H_{c}^{i}(X))$,
for each $i\geq0$, and this is a representation of $G$. The quantity

\[
\mathscr{L}(g,X):=\sum_{i\geq0}(-1)^{i}\Tr(g\mid H_{c}^{i}(X))=\Tr(g\mid H_{c}^{*}(X))\]
is called the \emph{Lefschetz number} of $X$ at $g$. A \emph{virtual
representation} of $G$ is an element in the Grothendieck group of
the semigroup generated by $\Irr(G)$ under the direct sum operation.
The function $\mathscr{L}(-,X):G\rightarrow\overline{\mathbb{Q}}_{l}$
is the character of the virtual representation $H_{c}^{*}(X)$ given
by the action of $G$ on $X$. Let $G$ be a finite group that acts
on the varieties $X$ and $Y$, respectively. Recall that we write
$X\sim Y$ if $H_{c}^{*}(X)=H_{c}^{*}(Y)$ as virtual $G$-representations.
We then have $X\sim Y$ if and only if $\mathscr{L}(-,X)=\mathscr{L}(-,Y)$,
and the relation $\sim$ is an equivalence relation. 
\begin{lem}
\label{lem:bijection}Suppose that $f:X\rightarrow Y$ is a (set-theoretic)
bijection between two varieties such that $f\phi=\phi f$, for some
Frobenius endomorphisms $\phi:X\rightarrow X$ and $\phi:Y\rightarrow Y$.
Let $g,g'$ be automorphisms of finite order of $X,Y$ such that $fg=g'f$.
Then $\mathscr{L}(g,X)=\mathscr{L}(g',Y)$.\end{lem}
\selectlanguage{british}%
\begin{proof}
As in the proof of \foreignlanguage{english}{\cite{dignemichel},
10.12 (ii), we have that for sufficiently large $m$,\[
|X^{g\phi^{m}}|=\sum_{y\in Y^{g'\phi^{m}}}|f^{-1}(y)^{g'\phi^{m}}|=|Y^{g'\phi^{m}}|,\]
which implies that $\mathscr{L}(g,X)=\mathscr{L}(g',Y)$.}
\end{proof}
Let $G$ be an affine algebraic group, and let $X\subseteq G$ be
a locally closed subset. Suppose that $H$ is a closed subgroup of
$G$, acting by multiplication on $G$, such that $X$ is stable under
the action of $H$. Then the quotient $X/H$ is a locally closed subset
of $G/H$. For a proof of this fact, see for example \cite{Springer-Bruhat-lemma},
Lemma~1.5. This shows that the quotient $X/H$ has a natural structure
of algebraic variety, which ensures that certain sets we will define
in the following are indeed varieties.

The following observations will be very useful in our analysis of
the cohomology of varieties. \foreignlanguage{english}{Let $G$ be
a finite group that acts on the variety $X$, and let $H\subset G$
be a subgroup such that there exists a $G$-equivariant morphism \[
\rho:X\longrightarrow G/H,\]
that is, $\rho$ satisfies $\rho(gx)=g\rho(x)$, for all $g\in G$,
$x\in X$. It then follows that $\rho$ is a surjection, and for any
$a\in G$, the stabiliser in $G$ of the fibre $\rho^{-1}(aH)$ is
$H\cap{}^{a}H$. Let $f$ be the fibre over the trivial coset $H\in G/H$.
Then every fibre of $\rho$ is isomorphic to $f$ via translation
by an element of $G$. Hence every $x\in X$ has the form $x=gy$,
for $g\in G$ and $y\in f$ which are uniquely determined up to the
action of $H$ given by $h(g,y)=(gh^{-1},hy)$. We thus have a $G$-equivariant
isomorphism\[
X\longiso(G\times f)/H,\qquad x\longmapsto(g,y)H.\]
Here $G$ acts on $(G\times f)/H$ via $g'(g,y)H=(g'g,y)H$. It follows
that \[
H_{c}^{*}((G\times f)/H)\cong\overline{\mathbb{Q}}_{l}[G]\otimes_{\overline{\mathbb{Q}}_{l}[H]}H_{c}^{*}(f)=\Ind_{H}^{G}H_{c}^{*}(f),\]
as virtual $G$-representations.}

\selectlanguage{english}%

\section{\label{sec:The-unramified-approach}The unramified approach}

Let $\mathbf{G}$ be a reductive group scheme over $\mathcal{O}_{F}$,
and let $r\geq1$ be an integer. A certain generalisation of the construction
of Deligne and Lusztig to the case $r\geq1$ was obtained by Lusztig
\cite{Lusztig-Fin-Rings} for $F$ of characteristic $p$, and in
\cite{Alex_Unramified_reps} for general $F$ and also for groups
over general finite local rings. The generalised Deligne-Lusztig varieties
in these constructions are attached to certain maximal tori in $\mathbf{G}\times\mathcal{O}_{F^{\mathrm{ur}}}$,
and are close analogues of the classical Deligne-Lusztig varieties.
Any maximal torus in $\mathbf{G}\times\mathcal{O}_{F^{\mathrm{ur}}}$
is an unramified torus in $\mathbf{G}\times_{\mathcal{O}_{F^{\mathrm{ur}}}}F^{\mathrm{ur}}$
in the sense that it splits over an unramified extension of $F$.
The construction given by these varieties can thus be seen as an {}``unramified''
generalisation of the construction of Deligne and Lusztig. We give
an outline of this construction.

Let $\phi:{G}_{r}\rightarrow{G}_{r}$ be a surjective endomorphism
of algebraic groups such that $G_{r}^{\phi}$ is finite. We call such
a map $\phi$ a \emph{Frobenius endomorphism}. Let $L:{G}_{r}\rightarrow{G}_{r}$,
denote the map $g\mapsto g^{-1}\varphi(g)$. Assume for simplicity
that $\mathbf{G}\times\mathcal{O}_{F^{\mathrm{ur}}}$  contains a
maximal torus $\mathbf{T}$ and a Borel subgroup $\mathbf{B}$ containing
$\mathbf{T}$, such that ${T}_{r}$ and ${B}_{r}$ are $\phi$-stable.
Let $\mathbf{U}$ be the unipotent radical of $\mathbf{B}$. By the
results in \cite{Alex-RedGreenAlg}, we know that\foreignlanguage{british}{
${B}_{r}$ is a self-normalising subgroup of ${G}_{r}$}. Note that
the assumption that ${B}_{r}$ be $\phi$-stable is not necessary
for the construction of the representations in \cite{Lusztig-Fin-Rings}
and \cite{Alex_Unramified_reps}, but it simplifies the models of
the varieties we consider here. 

\selectlanguage{british}%
Let $\mathcal{B}_{r}$ be the set of subgroups conjugate to ${B}_{r}$.
Since ${B}_{r}$ is self-normalising we have a bijection $\mathcal{B}_{r}\cong{G}_{r}/{B}_{r}$,
giving $\mathcal{B}_{r}$ a variety structure. As in the $r=1$ case,
we have a bijection\[
{G}_{r}\backslash(\mathcal{B}_{r}\times\mathcal{B}_{r})\longiso{B}_{r}\backslash{G}_{r}/{B}_{r}.\]
However, for $r>1$, the double ${B}_{r}$-${B}_{r}$ cosets are no
longer in one-to-one correspondence with elements of the group $N_{G_{r}}(T_{r})/T_{r}$,
and the structure of \foreignlanguage{english}{${B}_{r}\backslash{G}_{r}/{B}_{r}$
is too complex to admit any straightforward description.} Let $x\in{G}_{r}$
be an arbitrary element. In analogy with the $r=1$ case we can define
a variety\begin{multline*}
X_{r}(x):=\{B\in\mathcal{B}_{r}\mid(B,\phi(B))\in O(x)\}\\
\cong\{g\in{G}_{r}\mid g^{-1}\phi(g)\in{B}_{r}x{B}_{r}\}/{B}_{r},\\
\cong\{g\in{G}_{r}\mid g^{-1}\phi(g)\in x{B}_{r}\}/({B}_{r}\cap x{B}_{r}x^{-1}),\end{multline*}
where $O(x)$ denotes the orbit in $ $${G}_{r}\backslash(\mathcal{B}_{r}\times\mathcal{B}_{r})$
corresponding to the double coset ${B}_{r}x{B}_{r}$. In the same
way as for $r=1$, the finite group $G_{r}^{\phi}$ acts on $X_{r}(x)$
by left multiplication. For each \foreignlanguage{english}{$\hat{w}\in N_{G_{r}}(T_{r})$}
we also have a variety\begin{multline*}
\widetilde{X}_{r}(\hat{w}):=\{g\in{G}_{r}\mid g^{-1}\varphi(g)\in\hat{w}{U}_{r}\}/{U}_{r}\cap\hat{w}{U}_{r}\hat{w}^{-1}\\
=L^{-1}(\hat{w}{U}_{r})/{U}_{r}\cap\hat{w}{U}_{r}\hat{w}^{-1}.\end{multline*}
\foreignlanguage{english}{ The variety $\widetilde{X}_{r}(\hat{w})$
has a left action of ${G}_{r}^{\phi}$, and a commuting right action
of the group \[
{T}_{r}^{\hat{w}\phi}:=\{t\in{T}_{r}\mid\hat{w}\phi(t)\hat{w}^{-1}=t\}.\]
 }It is then not hard to verify, by the same method as for $r=1$,
that the varieties $\widetilde{X}_{r}(\hat{w})$ are finite $G_{r}^{\phi}$-covers
of $X_{r}(\hat{w})$. This depends on the fact that $\hat{w}$ normalises
the group ${T}_{r}$. The varieties $\widetilde{X}_{r}(\hat{w})$
(or rather, certain models isomorphic to them) were used in \cite{Lusztig-Fin-Rings}
and \cite{Alex_Unramified_reps} to construct certain generalised
Deligne-Lusztig representations. However, we will show in Subsection~\ref{sub:A-counterexample}
that the representations thus constructed leave out a non-trivial
subset of $\Irr(G_{r}^{\phi})$, for $r\geq2$. To remedy this situation
one would like to define further varieties that would produce the
missing representations. Given the above construction and the fact
that the elements \foreignlanguage{english}{$\hat{w}\in N_{G_{r}}(T_{r})$}
do not account for all of the double cosets in \foreignlanguage{english}{${B}_{r}\backslash{G}_{r}/{B}_{r}$},
it is a priori natural to define the following varieties \foreignlanguage{english}{(first
considered by Lusztig)}\[
L^{-1}(x{U}_{r})=\{g\in G_{r}\mid g^{-1}\phi(g)\in xU_{r}\},\quad\textrm{for any }x\in{G}_{r}.\]
Note that $L^{-1}(xU_{r})$ has an action of $U_{r}\cap xU_{r}x^{-1}$
by right multiplication, and the quotient $L^{-1}(xU_{r})/U_{r}\cap xU_{r}x^{-1}$
is a variety (see Section~\ref{sec:Preliminaries}). For $x=\hat{w}\in N_{{G}_{r}}({T}_{r})$
we have $L^{-1}(\hat{w}U_{r})/U_{r}\cap\hat{w}U_{r}\hat{w}^{-1}=\widetilde{X}_{r}(\hat{w})$,
and as we observed above, the variety $\widetilde{X}_{r}(\hat{w})$
is a finite cover of $X_{r}(\hat{w})$.\foreignlanguage{english}{
However, we point out that when }$x\notin N_{{G}_{r}}({T}_{r})$,
it is not in general the case that $L^{-1}(xU_{r})$, or even its
quotient $L^{-1}(xU_{r})/U_{r}\cap xU_{r}x^{-1}$, is a finite cover
of $X_{r}(x)$.\foreignlanguage{english}{ One might then hope that
in general any irreducible representation of ${G}_{r}^{\phi}$ is
realised by some variety $X_{r}(x)$ or $L^{-1}(xU_{r})$, for some
$x\in{G}_{r}$. This however, turns out to be not the case in general.
In the present section we will show that there exist irreducible representations
of $\SL_{2}(\mathcal{O}_{F,2})$, with $F$ of positive characteristic,
which are not realised in the cohomology of any variety of the form
$X_{2}(x)$ or $L^{-1}(xU_{2})$. Our proof proceeds as follows. First
we give an algebraic description of the irreducible representations
of $\SL_{2}(\mathcal{O}_{F,r})$, with particular emphasis on the
so-called nilpotent representations. We then analyse varieties of
the form $L^{-1}(xU_{2})$ and $X_{2}(x)$ and compare this to the
algebraic description of representations given earlier. Using computations
of Lusztig, giving the irreducible components of the cohomology of
$\widetilde{X}_{2}(\hat{w})$, where $B_{2}\hat{w}B_{2}\neq B_{2}$,
we can show that there exist representations in $\Irr(\SL_{2}(\mathcal{O}_{F,2}))$
which are not afforded by the varieties $L^{-1}(xU_{2})$ or $X_{2}(x)$.}

The following results will be applied in Subsection~\ref{sub:A-counterexample}
to the case where $\mathbf{G}=\SL_{2}$, $r=2$.\foreignlanguage{english}{ }
\selectlanguage{english}%
\begin{lem}
\label{lem:UxU}The inclusion $L^{-1}(x{U}_{r})\hookrightarrow L^{-1}({U}_{r}x{U}_{r})$
induces an isomorphism\[
L^{-1}(x{U}_{r})/{U}_{r}\cap x{U}_{r}x^{-1}\longiso L^{-1}({U}_{r}x{U}_{r})/{U}_{r},\]
commuting with the action of ${G}_{r}^{\phi}$ on both varieties.\end{lem}
\begin{proof}
Let $f$ be the composition of the maps \[
L^{-1}(x{U}_{r})\hookrightarrow L^{-1}({U}_{r}x{U}_{r})\rightarrow L^{-1}({U}_{r}x{U}_{r})/{U}_{r},\]
where the latter is the natural projection. Clearly $f$ is surjective,
because if $g{U}_{r}\in L^{-1}({U}_{r}x{U}_{r})/{U}_{r}$, with $L(g)\in uxu'$
for $u,u'\in{U}_{r}$, then $L(gu)=u^{-1}uxu'\varphi(u)\in x{U}_{r}$,
so $gu\in L^{-1}(x{U}_{r})$, and $f(gu)=g{U}_{r}$.

On the other hand, the fibre of $f$ at $g{U}_{r}$ is equal to\begin{multline*}
\{gv\in L^{-1}(x{U}_{r})\mid v\in{U}_{r}\}=\{gv\mid v^{-1}L(g)\varphi(v)\in x{U}_{r},\ v\in{U}_{r}\}\\
=\{gv\mid v^{-1}ux\in x{U}_{r},v\in{U}_{r}\}=\{gv\mid v^{-1}u\in{U}_{r}\cap x{U}_{r}x^{-1}\}\\
=\{gv\mid v=u\bmod{U}_{r}\cap x{U}_{r}x^{-1}\}.\end{multline*}
Factoring $L^{-1}(x{U}_{r})$ by ${U}_{r}\cap x{U}_{r}x^{-1}$ therefore
gives an isomorphism which commutes with the action of $G_{r}^{\phi}$.\end{proof}
\begin{lem}
\label{lem:F(x)UF(x)^-1}Let $x\in{G}_{r}$ be an arbitrary element,
and let $\lambda$ be an element such that $L(\lambda)=x$. Then there
is an isomorphism\[
L^{-1}(x{U}_{r})\longiso L^{-1}(\varphi(\lambda){U}_{r}\varphi(\lambda)^{-1}),\quad g\longmapsto g\lambda^{-1},\]
commuting with the action of $G_{r}^{\phi}$.\end{lem}
\begin{proof}
Let $g\in L^{-1}(x{U}_{r})$. Then \[
L(g\lambda^{-1})=\lambda L(g)\varphi(\lambda)^{-1}\in\lambda x{U}_{r}\varphi(\lambda)^{-1}=\varphi(\lambda){U}_{r}\varphi(\lambda)^{-1}.\]
It is clear that this map is a morphism of varieties, and it has an
obvious inverse.\end{proof}
\selectlanguage{british}%
\begin{lem}
\label{lem:xU sim wU}Suppose that $n\in N_{{G}_{r}}({T}_{r})$, and
let $x\in{B}_{r}n{B}_{r}$. Then \[
L^{-1}(x{U}_{r})/U_{r}\cap xU_{r}x^{-1}\sim L^{-1}(n{U}_{r})/U_{r}\cap nU_{r}n^{-1}.\]
\end{lem}
\selectlanguage{english}%
\begin{proof}
We can write $x$ as $utnt'u'$, for some $u,u'\in{U}_{r}$ and $t,t'\in{T}_{r}$.
Since $U_{r}$ is isomorphic to an affine space, \cite{dignemichel},
10.12 (ii) together with Lemma~\ref{lem:UxU} imply that\begin{multline*}
L^{-1}(x{U}_{r})/U_{r}\cap xU_{r}x^{-1}\sim L^{-1}({U}_{r}utnt'u'{U}_{r})\\
=L^{-1}({U}_{r}tnt'{U}_{r})\sim L^{-1}(tnt'{U}_{r})/U_{r}\cap tnt'U_{r}(tnt')^{-1}\\
=L^{-1}(t''n{U}_{r})/U_{r}\cap nU_{r}n^{-1},\end{multline*}
for some $t''\in{T}_{r}$. Since $t\mapsto n\varphi(t)n^{-1}$ is
a Frobenius map on ${T}_{r}$, The Lang-Steinberg theorem says that
there exists a $\lambda\in{T}_{r}$ such that $\lambda^{-1}n\varphi(\lambda)n^{-1}=t''$.
The map\begin{align*}
 & L^{-1}(t''n{U}_{r})/U_{r}\cap nU_{r}n^{-1}\longrightarrow L^{-1}(n{U}_{r})/U_{r}\cap nU_{r}n^{-1}\\
 & g(U_{r}\cap nU_{r}n^{-1})\longmapsto g\lambda^{-1}(U_{r}\cap nU_{r}n^{-1}),\end{align*}
is then an isomorphism of varieties which preserves the action of
${G}_{r}^{\phi}$. The lemma is proved. 
\end{proof}

\subsection{\label{sub:Reps-SL_2}The representations of $\SL_{2}(\mathcal{O}_{F,r})$}

Using results from Clifford theory and classification of conjugacy
orbits in certain algebras over the rings $\mathcal{O}_{F,r}$, it
is possible to completely describe the representations of the groups
$\SL_{2}(\mathcal{O}_{F,r})$, and $\GL_{2}(\mathcal{O}_{F,r})$.
In most cases, these algebras are the Lie algebras of the corresponding
group, with $\SL_{2}$, $p=2$ being a notable exception, as we will
see below. For $\SL_{2}$ with $p\neq2$ this method was employed
by Kutzko in his thesis (unpublished, see the announcement \cite{kutzkothesis})
and by Shalika (whose results remained unpublished until recently,
cf.~\cite{Shalika}). Around the same time the representations of
$\SL_{2}(\mathbb{Z}/p^{r}\mathbb{Z})$, including the case where $p=2$,
were also constructed by Nobs and Wolfart \cite{MR0444787,MR0444788},
by decomposing Weil representations. For $\GL_{2}$ with $\mathcal{O}_{F}=\mathbb{Z}_{p}$
and $p$ odd, the analogous result was given by Nagornyj \cite{nagornyj1},
and a general construction for all $\GL_{2}(\mathcal{O}_{F,r})$ can
be found in \cite{Alex_smooth_reps_GL2}. Recently, the $\SL_{2}$
case with $p\neq2$ was also reproduced in \cite{Jaikin-zeta}. We
will focus here on $\SL_{2}$, using the method of orbits and Clifford
theory, and without any restriction on $p$. The case where $p=2$
requires special treatment, and does not seem to have previously appeared
in the literature in this form. Proofs of the results we use can be
found in \cite{Shalika} and \cite{Alex_smooth_reps_GL2}, and we
will therefore omit details that can be found in these references.

Assume until the end of Subsection~\ref{sub:A-counterexample} that
$\mathbf{G}=\SL_{2}$, viewed as group scheme over $\mathcal{O}_{F}$.
\foreignlanguage{british}{Let $\mathbf{T}$ be the diagonal split
maximal torus in $\mathbf{G}$, $\mathbf{B}$ be the upper-triangular
Borel subgroup of $\mathbf{G}$, and $\mathbf{U}$ be the unipotent
radical of $\mathbf{B}$. Let $\mathbf{U}^{-}$ be the unipotent radical
of the Borel subgroup opposite to $\mathbf{B}$. }As usual, we identify
$G_{F,r}$ with the matrix group $\SL_{2}(\mathcal{O}_{F,r})$. Let
$\mathfrak{g}=\mathfrak{sl}_{2}$ be the Lie algebra of $\SL_{2}$,
viewed as a scheme over $\mathcal{O}_{F}$. Thus $\mathfrak{g}_{F,r}\cong\mathfrak{g}(\mathcal{O}_{F,r})$
is identified with the algebra of $2\times2$ matrices over $\mathcal{O}_{F,r}$
whose trace is zero. Assume first that $p\neq2$, and fix a natural
number $r>1$. For any natural number $i$ such that $r\geq i\geq1$
let $\rho_{r,i}:G_{F,r}\rightarrow G_{F,i}$ be the canonical surjective
homomorphism. For clarity, we will use the notation $K_{i}$ for the
kernel $G_{F,r}^{i}=\Ker\rho_{r,i}$. Assume from now on that $i\geq r/2$.
Then $K_{i}=1+\mathfrak{p}_{F}^{i}\mathfrak{g}_{F,r-i}$ and the map
$x\mapsto1+\varpi^{i}x$ induces an isomorphism $\mathfrak{g}_{F,r-i}\iso K_{i}$.
The group $G_{F,r}$ acts on $\mathfrak{g}_{F,r-i}$ by conjugation,
via its quotient $G_{F,r-i}$. This action is transformed by the above
isomorphism into the action of $G_{F,r}$ on the normal subgroup $K_{i}$. 

Fix an additive character $\psi:\mathcal{O}_{F}\rightarrow\overline{\mathbb{Q}}_{l}\vphantom{\overline{\mathbb{Q}}}^{\hspace{-3pt}\times}$
with conductor $\mathfrak{p}_{F}^{r}$, and define for any $\beta\in\mathfrak{g}_{F,r-i}$
a character $\psi_{\beta}:K_{i}\rightarrow\overline{\mathbb{Q}}_{l}\vphantom{\overline{\mathbb{Q}}}^{\hspace{-3pt}\times}$
by\[
\psi_{\beta}(x)=\psi(\Tr(\beta(x-1))).\]
Then $\beta\mapsto\psi_{\beta}$ gives an isomorphism \[
\mathfrak{g}_{F,r-i}\cong\textrm{Hom}(K_{i},\overline{\mathbb{Q}}_{l}\vphantom{\overline{\mathbb{Q}}}^{\hspace{-3pt}\times}),\]
and for $g\in G_{F,r}$, we have $\rho_{r-i}(g)\beta\rho_{r-i}(g)^{-1}\mapsto(\psi_{\beta})^{g}$.

Set $l=[\frac{r+1}{2}]$, $l'=[\frac{r}{2}]$; thus $l+l'=r$. Let
$\pi$ be an irreducible representation of $G_{F,r}$. By Clifford's
theorem, restricting $\pi$ to $K_{l}$ determines an orbit of characters
on $K_{l}$, and hence (by the above isomorphism) an orbit in $\mathfrak{g}_{F,l'}$.
If the orbit is in $\mathfrak{p}_{F}\mathfrak{g}_{F,l'}$, then $\pi$
is trivial on $K_{r-1}$, and so factors though $G_{F,r-1}$. We are
only concerned with \emph{primitive} representations, that is, those
which do not factor through $G_{F,r-1}$. It is therefore enough to
consider orbits in $\mathfrak{g}_{F,l'}\setminus\mathfrak{p}_{F}\mathfrak{g}_{F,l'}$.
For any natural number $r'$ such that $r\geq r'\geq1$ we call an
element $\beta\in\mathfrak{g}_{F,r'}$ \emph{regular} if the centraliser
$C_{G_{1}}(\rho_{r',1}(\beta))$ in $G_{1}\cong\mathbf{G}(\overline{\mathbb{F}}_{q})$
is abelian. We then have $C_{G_{r'}}(\beta)=\mathcal{O}_{r'}[\beta]\cap G_{r'}$,
in the connected algebraic group $G_{r'}$. The orbits in $\mathfrak{g}_{F,l'}\setminus\mathfrak{p}_{F}\mathfrak{g}_{F,l'}$
can be easily classified thanks to the fact that they are all regular.
More precisely, the orbits in $\mathfrak{g}_{F,l'}\setminus\mathfrak{p}_{F}\mathfrak{g}_{F,l'}$
are of three basic types, according to their reductions mod $\mathfrak{p}_{F}$:
There are the orbits with split characteristic polynomial and distinct
eigenvalues mod $\mathfrak{p}_{F}$, the ones which have irreducible
characteristic polynomial mod $\mathfrak{p}_{F}$, and those which
are nilpotent mod $\mathfrak{p}_{F}$. The primitive representations
of these three types are called \emph{split, cuspidal, }and \emph{nilpotent,}
respectively.

The construction of the representations of $G_{F,r}$ with a given
orbit $\Omega\in\mathfrak{g}_{F,l'}\setminus\mathfrak{p}_{F}\mathfrak{g}_{F,l'}$
proceeds as follows. Pick a representative $\beta\in\Omega$, and
consider the corresponding character $\psi_{\beta}$ on $K_{l}$.
The stabiliser in $G_{F,r}$ of $\psi_{\beta}$ is given by\[
\Stab_{G_{F,r}}(\psi_{\beta})=C_{G_{F,r}}(\hat{\beta})K_{l'},\]
where $\hat{\beta}\in\mathfrak{g}_{F,r}$ is an element such that
$\rho_{r,l}(\hat{\beta})=\beta$. Assume first that $r$ is even so
that $l=l'$. Since $C_{G_{F,r}}(\hat{\beta})$ is abelian, the character
$\psi_{\beta}$ can be extended to a character on $\Stab_{G_{F,r}}(\psi_{\beta})$,
and all the irreducible representations of $\Stab_{G_{F,r}}(\psi_{\beta})$
containing $\psi_{\beta}$ are obtained in this way. Inducing a representation
of $\Stab_{G_{F,r}}(\psi_{\beta})$ containing $\psi_{\beta}$ to
$G_{F,r}$ gives an irreducible representation, and it is clear that
we get all the irreducible representations of $G_{F,r}$ with orbit
$\Omega$ in this way.

Now assume that $r$ is odd. In this case there are several equivalent
variations of the construction, but they all involve (at least for
some orbits) a step where a representation of a group is shown to
have a unique representation lying above it in a larger group. The
other steps consist of various lifts and induction from $\Stab_{G_{F,r}}(\psi_{\beta})$,
as in the case for $r$ even. For full details, see \cite{Shalika}
for $\SL_{2}$, and \cite{Alex_smooth_reps_GL2} for the closely related
case of $\GL_{2}$, respectively.

Now consider the case where $p=2$. In this case the association $\beta\mapsto\psi_{\beta}$
does no longer give an isomorphism between $\mathfrak{g}_{F,r-i}$
and the character group of $K_{i}$. To remedy this, we first consider
the analogous situation for $\GL_{2}$ where the role of $\mathfrak{g}_{F,r-i}$
is played by the matrix algebra $\M_{2}(\mathcal{O}_{F,r-i})$, and
the analogous map $\beta\mapsto\psi_{\beta}$ is indeed an isomorphism
(for any $p$). The $i^{\text{th}}$ congruence kernel in $\GL_{2}(\mathcal{O}_{F,r})$
has the form $1+\mathfrak{p}_{F}^{i}\M_{2}(\mathcal{O}_{F,r-i})$,
and so it contains $K_{i}$ as a subgroup of index $|\mathcal{O}_{F,r-i}|$.
For every $\beta\in\M_{2}(\mathcal{O}_{F,r-i})$ we have a character
$\psi_{\beta}|_{K_{i}}$ obtained by restricting the character $\psi_{\beta}$
on $1+\mathfrak{p}_{F}^{i}\M_{2}(\mathcal{O}_{F,r-i})$ to $K_{i}$.
Then $\beta\mapsto\psi_{\beta}|_{K_{i}}$ is obviously a surjective
homomorphism $\M_{2}(\mathcal{O}_{F,r-i})\rightarrow\textrm{Hom}(K_{i},\overline{\mathbb{Q}}_{l}\vphantom{\overline{\mathbb{Q}}}^{\hspace{-3pt}\times})$.
It is easily seen that the kernel of this homomorphism is the subgroup
$Z_{r-1}$ of scalar matrices in $\M_{2}(\mathcal{O}_{F,r-i})$. We
therefore have an isomorphism\[
\M_{2}(\mathcal{O}_{F,r-i})/Z_{r-i}\longiso\textrm{Hom}(K_{i},\overline{\mathbb{Q}}_{l}\vphantom{\overline{\mathbb{Q}}}^{\hspace{-3pt}\times}),\quad\beta+Z_{r-i}\mapsto\psi_{\beta}|_{K_{i}}.\]
Since $Z_{r-i}$ is centralised by $G_{F,r}$, we see that for any
$g\in G_{F,r}$, we have \[
\rho_{r-i}(g)\beta\rho_{r-i}(g)^{-1}\mapsto(\psi_{\beta}|_{K_{i}})^{g}.\]
 As before, let $l=[\frac{r+1}{2}]$, $l'=[\frac{r}{2}]$. If $\beta\in\mathfrak{p}_{F}\M_{2}(\mathcal{O}_{F,l'})/Z_{l'}$,
then $\psi_{\beta}|_{K_{l}}$ is trivial on $K_{r-1}$, and so an
irreducible representation of $G_{F,r}$ whose restriction to $K_{l}$
contains this $\psi_{\beta}|_{K_{l}}$ must factor through $G_{F,r-1}$,
and hence is not primitive. To construct the primitive representations,
the first task is now to classify the orbits under the action of $G_{F,r}$
on $\M_{2}(\mathcal{O}_{F,l'})/Z_{l'}\setminus\mathfrak{p}_{F}\M_{2}(\mathcal{O}_{F,l'})/Z_{l'}$.
The following is a list a representatives of these orbits:
\selectlanguage{british}%
\begin{enumerate}
\item $\begin{pmatrix}a & 0\\
0 & 0\end{pmatrix}$, $a\in\mathcal{O}_{F,l'}^{\times}$,
\item $\begin{pmatrix}0 & 1\\
\Delta & s\end{pmatrix}$, where $\Delta,s\in\mathcal{O}_{F,l'}$, and $x^{2}-sx-\Delta$ is
irreducible mod $\mathfrak{p}_{F}$,
\item $\begin{pmatrix}0 & 1\\
\Delta & s\end{pmatrix}$, where $\Delta,s\in\mathfrak{p}_{F}$.
\end{enumerate}
The construction of representations then proceeds as in the case $p\neq2$. 
\begin{rem*}
Clearly the method used in the case $p=2$ could also be applied when
$p\neq2$. We have however chosen to give the two separate cases in
order to illustrate their contrasts. Note that when $p\neq2$ the
embedding $\mathfrak{g}_{F,l'}\hookrightarrow\M_{2}(\mathcal{O}_{F,l'})$
induces a $G_{F,r}$-equivariant isomorphism \[
\mathfrak{g}_{F,l'}\longiso\M_{2}(\mathcal{O}_{F,l'})/Z_{l'},\]
so in general the algebra $\M_{2}(\mathcal{O}_{F,l'})/Z_{l'}$ is
the right object, rather than the Lie algebra $\mathfrak{g}_{F,l'}$,
in which to consider orbits.
\end{rem*}
In the following we will be especially interested in the \emph{nilpotent}
representations of $G_{F,2}\cong\SL_{2}(\mathcal{O}_{F,2})$, that
is, the irreducible primitive representations whose orbits mod $\mathfrak{p}_{F}$
are nilpotent, or contain a nilpotent element mod $Z_{1}$ when $p=2$,
respectively. We call the corresponding orbits nilpotent (although
in the $p=2$ case, they are strictly speaking only nilpotent mod
centre). The construction of representations given above shows that
the nilpotent representations are induced from 1-dimensional representations
on $\Stab_{G_{F,2}}(\psi_{\beta}|_{K_{1}})$, where $\beta$ is a
representative of a nilpotent orbit. When $p\neq2$ there are \foreignlanguage{english}{exactly
two nilpotent orbits} in $ $\foreignlanguage{english}{$\mathfrak{g}_{F,1}\setminus\mathfrak{p}_{F}\mathfrak{g}_{F,1}$,
given by the representatives\[
\begin{pmatrix}0 & 1\\
0 & 0\end{pmatrix},\quad\begin{pmatrix}0 & \zeta\\
0 & 0\end{pmatrix},\]
respectively (here $\zeta\in\mathbb{F}_{q}^{\times}$ is a non-square
element). When $p=2$ there is just one nilpotent-mod-$Z_{1}$ orbit
in $\M_{2}(\mathcal{O}_{F,1})/Z_{1}\setminus\mathfrak{p}_{F}\M_{2}(\mathcal{O}_{F,1})/Z_{1}$},
given by the representative $\left(\begin{smallmatrix}0 & 1\\
0 & 0\end{smallmatrix}\right)$. If we let $\beta$ be any of these representatives, then the stabiliser
$\Stab_{G_{F,2}}(\psi_{\beta}|_{K_{1}})$ is given by\foreignlanguage{english}{
}\[
S:=\Stab_{G_{F,2}}(\psi_{\beta}|_{K_{1}})=\{\pm1\}U_{F,2}K_{1},\]
where $\{\pm1\}$ denotes a subgroup of scalar matrices (which is
equal to the centre of $G_{F,2}$ for $p\neq2$, and is trivial for
$p=2$), and $U_{F,2}$ is isomorphic to the subgroup of $\mathbf{G}(\mathcal{O}_{F,2})$
of upper unitriangular matrices. The index of $S$ in $G_{F,2}$ is
equal to $(q^{2}-1)/2$ when $p\neq2$, and equal to $q^{2}-1$ when
$p=2$. It is not hard to show that the commutator subgroup of $S$
is $[S,S]=B_{F,2}^{1}=B_{F,2}\cap K_{1}$. Thus all nilpotent representations
of $G_{F,2}$ are components of the induced representation $\Ind_{B_{F,2}^{1}}^{G_{F,2}}\mathbf{1}$.
Each $\psi_{\beta}$ has $|S/K_{1}|$ extensions to $S$, and each
such extension induces to a distinct nilpotent representation. When
$p\neq2$ we thus have $4q$ nilpotent representations, all of which
have dimension $(q^{2}-1)/2$. When $p=2$ we have $q$ nilpotent
representations, all of which have dimension $q^{2}-1$.

We will have occasion to consider the question of which nilpotent
representations occur as components of $\Ind_{U_{F,2}}^{G_{F,2}}\mathbf{1}$.
By the above we know that any nilpotent representation of $G_{F,2}$
is of the form $\Ind_{S}^{G_{F,2}}\rho$, for some $\rho$ such that
$\rho|K_{1}$ contains $\psi_{\beta}$, with $\beta$ one of the above
nilpotent representatives. By Mackey's intertwining number formula,
we have\[
\langle\Ind_{S}^{G_{F,2}}\rho,\Ind_{U_{F,2}}^{G_{F,2}}\mathbf{1}\rangle=\sum_{x\in S\backslash G_{F,2}/U_{F,2}}\langle\rho|_{S\cap\vphantom{U_{F,2}}^{x}U_{F,2}},\mathbf{1}\rangle,\]
and since $S$ contains $K_{1}$ we can identify $S\backslash G_{F,2}/U_{F,2}$
with $U_{F,1}\backslash G_{F,1}/U_{F,1}$. To calculate the value
of the right-hand side it is thus enough to let $x$ run through elements
in $T_{F,2}$ and elements in $\hat{w}T_{F,2}$, respectively ($\hat{w}\in N_{G_{F,2}}(T_{F,2})$
denotes a lift of the non-trivial element of the Weyl group of $\SL_{2}(k)$).
Since $T_{F,2}$ normalises $U_{F,2}$, it is moreover enough to consider
only $x=1$ and $x=\hat{w}$. For $x=1$ we get a term $\langle\rho|_{U_{F,2}},\mathbf{1}\rangle$,
and for $x=\hat{w}$ we get a term $\langle\rho|_{(U^{-})_{F,2}^{1}},\mathbf{1}\rangle$.
The latter is always zero, since $\rho|_{(U^{-})_{F,2}^{1}}=\psi_{\beta}|_{(U^{-})_{F,2}^{1}}\neq\mathbf{1}$
for our choice of $\beta$. Hence we conclude that $\Ind_{S}^{G_{F,2}}\rho$
is contained in $\Ind_{U_{F,2}}^{G_{F,2}}\mathbf{1}$ if and only
$\langle\rho|_{U_{F,2}},\mathbf{1}\rangle=1$. In particular, since
there exist representations of $S$ which are lifts of $\psi_{\beta}$
and which are non-trivial on $U_{F,2}$, we see that there exist nilpotent
representations which are not components of $\Ind_{U_{F,2}}^{G_{F,2}}\mathbf{1}$.

\selectlanguage{english}%

\subsection{\label{sub:A-counterexample}Inadequacy of the unramified varieties}

We keep the assumption $\mathbf{G}=\SL_{2}$ until the end of this
subsection. We will show that there exist nilpotent representations
of $G_{F,2}$ which cannot be realised as components of the cohomology
of varieties of the form $L^{-1}(xU_{2})$, $L^{-1}(xU_{2})/U_{2}\cap xU_{2}x^{-1}$,
or $X_{2}(x)$, for $x\in{G}_{2}$. More precisely, we show that the
only nilpotent representations which can be realised in this way are
the irreducible components of $\Ind_{U_{F,2}}^{G_{F,2}}\mathbf{1}$.
As we saw above, these do not account for all the nilpotent representations
of $G_{F,2}$.

Let $\phi:G_{2}\rightarrow G_{2}$ be the standard Frobenius endomorphism
induced by the map which sends every matrix entry to its $q^{\mathrm{th}}$
power. Then $G_{F,2}=G_{2}^{\phi}$, and we will use either of these
ways of writing the group, depending on the context. Moreover, each
of the subgroups $T_{2}$, $B_{2}$, $U_{2}$, and $(U^{-})_{2}$
is $\phi$-stable. We need a description of the double cosets ${B}_{2}\backslash{G}_{2}/{B}_{2}$.
One checks directly that a set of representatives is given by\[
\Bigl{\{}1,\ w:=\begin{pmatrix}0 & 1\\
-1 & 0\end{pmatrix},\ e:=\begin{pmatrix}1 & 0\\
\varpi & 1\end{pmatrix}\Bigr\}.\]
Note that $e\in(U^{-})_{2}^{1}$ and that for any $a\in(U^{-})_{2}^{1}-\{1\}$,
we have $U_{2}\cap aU_{2}a^{-1}=U_{2}^{1}$, which is an affine space.
In this case, \cite{dignemichel}, 10.12 (ii) implies that $L^{-1}(aU_{2})\sim L^{-1}(aU_{2})/U_{2}^{1}$.
Note also that $U_{2}\cap wU_{2}w^{-1}=\{1\}$.
\begin{prop}
\label{pro:{1,w}(U^-)}Let $x\in G_{2}$ be an arbitrary element.
Then there exists an element $y\in\{1,w\}\cup(U^{-})_{2}^{1}$ such
that $L^{-1}(xU_{2})\sim L^{-1}(yU_{2})$.\end{prop}
\begin{proof}
The elements $1$ and $w$ normalise ${T}_{2}$ so, by Lemma~\ref{lem:xU sim wU},
for any element $x\in{B}_{2}$ we have $L^{-1}(xU_{2})\sim L^{-1}(U_{2})$,
and for any $x\in{B}_{2}w{B}_{2}$ we have $L^{-1}(xU_{2})\sim L^{-1}(wU_{2})$.

In contrast, no element in $B_{2}eB_{2}$ normalises ${T}_{2}$. Assume
that $x=utet'u'$, where $u,u'\in{U}_{2}$ and $t,t'\in{T}_{2}$.
Then $L^{-1}(utet'u'{U}_{2})\sim L^{-1}({U}_{2}tet'{U}_{2})\sim L^{-1}(tet'{U}_{2})$,
and by Lemma \ref{lem:F(x)UF(x)^-1} we have $L^{-1}(tet'{U}_{2})\sim L^{-1}(\varphi(\lambda){U}_{2}\varphi(\lambda)^{-1})$,
where $\lambda\in{G}_{2}$ is such that $L(\lambda)=tet'$. Since
$tet'\in(U^{-})_{2}^{1}T_{2}$ and the group $(U^{-})_{2}^{1}$ is
$\phi$-stable, we can take $\lambda\in(U^{-})_{2}^{1}T_{2}$, by
the Lang-Steinberg theorem. Writing $\lambda=vs$, with some $v\in(U^{-})_{2}^{1}$
and $s\in{T}_{2}$, we get \begin{multline*}
L^{-1}(\varphi(\lambda){U}_{2}\varphi(\lambda)^{-1})=\\
L^{-1}(\varphi(vs){U}_{2}\varphi(vs)^{-1})=L^{-1}(\varphi(v){U}_{2}\varphi(v)^{-1})\sim L^{-1}(L(v){U}_{2}).\end{multline*}
Since the group $(U^{-})_{2}^{1}$ is $\phi$-stable, we have $L(v)=v{}^{-1}\varphi(v)\in(U^{-})_{2}^{1}$.
Hence, for every $x\in{B}_{2}e{B}_{2}$, we have $L^{-1}(xU_{2})\sim L^{-1}(y{U}_{2})$,
for some $y\in(U^{-})_{2}^{1}$. \end{proof}
\begin{thm}
\label{thm:Inadequacy}Let $y\in(U^{-})_{2}^{1}$\foreignlanguage{british}{.}
Then $L^{-1}(yU_{2})\sim\widetilde{X}_{2}(1)$, and hence \[
H_{c}^{*}(L^{-1}(yU_{2}))\cong\Ind_{U_{2}^{\phi}}^{G_{2}^{\phi}}\mathbf{1}\]
as $G_{2}^{\phi}$-representations.\end{thm}
\begin{proof}
We use the observations from the end of Section~\ref{sec:Preliminaries}.
Consider the composition of the maps\[
\rho:L^{-1}(yU_{2})/{U}_{2}^{1}\xrightarrow{\ \rho_{2,1}\ }X_{1}({U}_{1})\longrightarrow G_{1}^{\phi}/U_{1}^{\phi}\cong G_{2}^{\phi}/U_{2}^{\phi}(G_{2}^{1})^{\phi},\]
where the first map is the restriction of $\rho_{2,1}:G_{2}\rightarrow G_{1}$,
and the second map is given by $g\mapsto gU_{1}^{\phi}$. Then $\rho$
is clearly $G_{2}^{\phi}$-equivariant. The fibre $f:=\rho^{-1}(U_{2}^{\phi}(G_{2}^{1})^{\phi})$
over the trivial coset in $G_{2}^{\phi}/U_{2}^{\phi}(G_{2}^{1})^{\phi}$
is given by\[
f=\{um\in{U}_{2}{G}_{2}^{1}\mid(um)^{-1}\phi(um)\in y{U}_{2}\}/{U}_{2}^{1}.\]
Pick a $\lambda\in(U^{-})_{2}^{1}$ such that $\lambda^{-1}\phi(\lambda)=y$.
Then the translation $x\mapsto x\lambda^{-1}$ induces a $U_{2}^{\phi}(G_{2}^{1})^{\phi}$-equivariant
isomorphism\[
f\longiso f\lambda^{-1}=\{um\in{U}_{2}{G}_{2}^{1}\mid(um)^{-1}\phi(um)\in\phi(\lambda){U}_{2}\phi(\lambda)^{-1}\}/{U}_{2}^{1}.\]
We now observe that the group $\phi(\lambda){U}_{2}\phi(\lambda)^{-1}$
is contained in ${U}_{2}{T}_{2}^{1}$. Thus, every element in $f\lambda^{-1}$
is $\phi$-fixed up to right multiplication by some element in ${U}_{2}{T}_{2}^{1}$.
Hence there is a map\[
\rho':f\lambda^{-1}\longrightarrow({U}_{2}{G}_{2}^{1}/{U}_{2}{T}_{2}^{1})^{\phi}\cong U_{2}^{\phi}(G_{2}^{1})^{\phi}/U_{2}^{\phi}(T_{2}^{1})^{\phi},\qquad x\longmapsto xU_{2}^{\phi}(G_{2}^{1})^{\phi},\]
which is clearly $U_{2}^{\phi}(G_{2}^{1})^{\phi}$-equivariant. Define
$f'$ to be the fibre of $\rho'$ over the trivial coset. Then\[
f'=\{um\in{U}_{2}{T}_{2}^{1}\mid(um)^{-1}\phi(um)\in\phi(\lambda){U}_{2}\phi(\lambda)^{-1}\}/{U}_{2}^{1},\]
which has a left action of $U_{2}^{\phi}(T_{2}^{1})^{\phi}$, and
a right action of $(T_{2}^{1})^{\phi}$.

\selectlanguage{british}%
We now show that the \foreignlanguage{english}{$U_{2}^{\phi}(T_{2}^{1})^{\phi}$}-representation
afforded by $f'$ is isomorphic to $\Ind_{U_{2}^{\phi}}^{U_{2}^{\phi}(T_{2}^{1})^{\phi}}\mathbf{1}$.
Define the variety\[
V=\{g\in{U}_{2}{T}_{2}^{1}\mid g^{-1}\phi(g)\in{U}_{2}\}={U}_{2}({T}_{2}^{1})^{\phi}.\]
This has a left action of $U_{2}^{\phi}(T_{2}^{1})^{\phi}$ and a
right action of ${U}_{2}^{\phi}$. We have $V/{U}_{2}\cong U_{2}^{\phi}(T_{2}^{1})^{\phi}/{U}_{2}^{\phi}$,
so $V$ affords the representation $\Ind_{U_{2}^{\phi}}^{U_{2}^{\phi}(T_{2}^{1})^{\phi}}\mathbf{1}$,
that is \[
H_{c}^{*}(V)\cong\Ind_{U_{2}^{\phi}}^{U_{2}^{\phi}(T_{2}^{1})^{\phi}}\mathbf{1},\]
as $U_{2}^{\phi}(T_{2}^{1})^{\phi}$-representations. Now, for every
$u\in{U}_{2}$ there exists a $t_{u}\in{T}_{2}^{1}$ such that $ut_{u}\in f'$,
and this $t_{u}$ is unique up to multiplication by $({T}_{2}^{1})^{\phi}$.
Hence, by choosing such a $t_{um}$ for each $um\in f'$, we can write
each element in $f'$ uniquely in the form $ut_{u}a$, where $u\in{U}_{2}^{1}$,
$t_{u}\in{T}_{2}^{1}$, and $a\in({T}_{2}^{1})^{\phi}$. Moreover,
we may always choose the same $t_{u}$ for all elements $vsus^{-1}$,
where $v\in{U}_{2}^{\phi}$ and $s\in({T}_{2}^{1})^{\phi}$. Similarly,
we may always choose $t_{u}$ so that $\phi^{m}(t_{u})=t_{\phi^{m}(u)}$,
for all natural numbers $m\geq1$. We can then define a bijective
function \[
\eta:f'\longrightarrow V,\qquad ut_{u}a\longmapsto ua.\]
For $vs\in U_{2}^{\phi}(T_{2}^{1})^{\phi}$ we have \[
\eta(vsut_{u}a)=\eta(v(sus^{-1}t_{u}sa))=v(sus^{-1})sa=vsua,\]
so $\eta$ is $U_{2}^{\phi}(T_{2}^{1})^{\phi}$-equivariant. Let $m$
be a natural number such that $\phi^{m}(\lambda)=\lambda$. Then $\phi^{m}$
is a Frobenius endomorphism on $f'$. Furthermore, $\phi^{m}$ is
clearly a Frobenius endomorphism which stabilises $V$. The bijection
$\eta$ satisfies \[
\eta(\phi^{m}(ut_{u}a))=\eta(\phi^{m}(u)\phi^{m}(t_{u})a)=\eta(\phi^{m}(u)t_{\phi^{m}(u)}a)=\phi^{m}(u)a=\phi^{m}(ua),\]
so $\eta$ commutes with the Frobenius endomorphisms $\phi^{m}$ on
$f'$ and $V$, respectively. By Lemma~\ref{lem:bijection} $f'$
and $V$ afford the same $U_{2}^{\phi}(T_{2}^{1})^{\phi}$-representation,
and so\begin{multline*}
H_{c}^{*}(L^{-1}(yU_{2}))\\
\cong\Ind_{U_{2}^{\phi}(G_{2}^{1})^{\phi}}^{G_{2}^{\phi}}\Ind_{U_{2}^{\phi}(T_{2}^{1})^{\phi}}^{U_{2}^{\phi}(G_{2}^{1})^{\phi}}\Ind_{U_{2}^{\phi}}^{U_{2}^{\phi}(T_{2}^{1})^{\phi}}\mathbf{1}=\Ind_{U_{2}^{\phi}}^{G_{2}^{\phi}}\mathbf{1}\\
\cong H_{c}^{*}(\widetilde{X}_{2}(1)).\end{multline*}

\end{proof}
The representations realised by the variety $\widetilde{X}_{2}(1)$,
that is, the irreducible components of \foreignlanguage{british}{$\Ind_{U_{2}^{\phi}}^{G_{2}^{\phi}}\mathbf{1}$,}
are just the irreducible components of the representations obtained
by lifting characters of $T_{2}^{\phi}$ to $B_{2}^{\phi}$, and inducing
to $G_{2}^{\phi}$. As we saw in the end of Section \ref{sub:Reps-SL_2},
not all of the nilpotent representations are of this form.

When $F$ is a local field of characteristic $p$, Lusztig \cite{Lusztig-Fin-Rings}
has identified the representations realised by the variety $\widetilde{X}_{2}(w)$.
In particular, none of them is of dimension $(q^{2}-1)/2$ when $p\neq2$,
or of dimension $q^{2}-1$ when $p=2$, so in this case the variety
$\widetilde{X}_{2}(w)$ does not realise any of the nilpotent representations
of $G_{2}^{\phi}=G_{F,2}$. Thus the results of this section imply
that there are nilpotent representations of $\mbox{SL}_{2}(\mathbb{F}_{q}[[\varpi]]/(\varpi^{2}))$
which are not realised in the cohomology of any of the varieties $L^{-1}(xU_{2})$,
or equivalently, the varieties $L^{-1}(xU_{2})/U_{2}\cap xU_{2}x^{-1}$,
for $x\in{G}_{2}$.
\begin{rem*}
It seems likely that Lusztig's result on the representations afforded
by $\widetilde{X}_{2}(w)$ hold in any characteristic, in particular,
that $\widetilde{X}_{2}(w)$ does not afford any nilpotent representation
of $G_{F,2}$, for any non-archimedean local field $F$. More precisely,
every irreducible representation of $G_{F,2}$ afforded by $\widetilde{X}_{2}(w)$
should be either non-primitive or cuspidal. Since the results in this
section hold uniformly in any characteristic, this would imply the
inadequacy of the varieties $L^{-1}(xU_{2})$ also for the group $\SL_{2}(\mathbb{Z}/p^{r}\mathbb{Z})$. 
\end{rem*}
As we remarked in the beginning of the section, the variety $L^{-1}(eU_{2})/U_{2}^{1}$
is not a finite cover of $X_{2}(e)$, so the representations afforded
by the latter are not necessarily all afforded by the former (as is
the case for the covers $\widetilde{X}_{r}(\hat{w})$ of $X_{r}(\hat{w})$,
for $\hat{w}\in N_{G_{r}}(T_{r})$). It is thus a priori conceivable
that $X_{2}(e)$ may yield further representations not obtainable
by $L^{-1}(eU_{2})$. The following result shows that this is not
the case.
\begin{prop}
\label{pro:Inadeq-X_2(e)}We have\[
H_{c}^{*}(X_{2}(e))=\Big(\Ind_{B_{2}^{\phi}(G_{2}^{1})^{\phi}}^{G_{2}^{\phi}}\mathbf{1}\Big)-\Ind_{B_{2}^{\phi}}^{G_{2}^{\phi}}\mathbf{1},\]
as virtual $G_{2}^{\phi}$-representations.\end{prop}
\begin{proof}
Consider the composition of the maps\[
X_{2}(e)\xrightarrow{\ \rho_{2,1}\ }L^{-1}({B}_{1})/{B}_{1}\longiso G_{1}^{\phi}/B_{1}^{\phi}\longiso G_{2}^{\phi}/B_{2}^{\phi}(G_{2}^{1})^{\phi}.\]
This gives a $G_{2}^{\phi}$-equivariant map $X_{2}(e)\rightarrow G_{2}^{\phi}/B_{2}^{\phi}(G_{2}^{1})^{\phi}$.
The fibre of the trivial coset under this map is\[
f:=\{g\in{B}_{2}{G}_{2}^{1}\mid g^{-1}\phi(g)\in{B}_{2}e{B}_{2}\}/{B}_{2}.\]
Thus we have \[
H_{c}^{*}(X_{2}(e))=\Ind_{B_{2}^{\phi}(G_{2}^{1})^{\phi}}^{G_{2}^{\phi}}H_{c}^{*}(f).\]
Now an element in ${B}_{2}^{\phi}({G}_{2}^{1})^{\phi}$ must lie in
exactly one of the double cosets ${B}_{2}$ and ${B}_{2}e{B}_{2}$.
Hence \[
f\sqcup\{g\in{B}_{2}{G}_{2}^{1}\mid g^{-1}\phi(g)\in{B}_{2}\}/{B}_{2}={B}_{2}{G}_{2}^{1}/{B}_{2}.\]
Since ${B}_{2}{G}_{2}^{1}/{B}_{2}\cong{G}_{2}^{1}/{B}_{2}^{1}$ is
an affine space, the $G_{2}^{\phi}$-representation afforded by it
is the trivial representation. Moreover, the variety \[
\{g\in{B}_{2}{G}_{2}^{1}\mid g^{-1}\phi(g)\in{B}_{2}\}/{B}_{2}\]
 is isomorphic to ${B}_{2}^{\phi}({G}_{2}^{1})^{\phi}/{B}_{2}^{\phi}$,
and so affords the representation $\Ind_{B_{2}^{\phi}}^{B_{2}^{\phi}(G_{2}^{1})^{\phi}}\mathbf{1}$.
Putting these results together, we get\[
H_{c}^{*}(f\sqcup\{g\in{B}_{2}{G}_{2}^{1}\mid g^{-1}\phi(g)\in{B}_{2}\}/{B}_{2})=H_{c}^{*}(f)+\Ind_{B_{2}^{\phi}}^{B_{2}^{\phi}(G_{2}^{1})^{\phi}}\mathbf{1}=\mathbf{1},\]
whence the result.
\end{proof}
The irreducible components of the representation $\Ind_{B_{2}^{\phi}(G_{2}^{1})^{\phi}}^{G_{2}^{\phi}}\mathbf{1}$
are all non-primi\-tive, since they have $(G_{2}^{1})^{\phi}$ in
their respective kernels. Moreover, the irreducible components of
$\Ind_{B_{2}^{\phi}}^{G_{2}^{\phi}}\mathbf{1}$ form a subset of the
irreducible components of $\Ind_{U_{2}^{\phi}}^{G_{2}^{\phi}}\mathbf{1}$.
Thus, the variety $X_{2}(e)$ does not afford any nilpotent representations
of $G_{2}^{\phi}=G_{F,2}$ which are not already afforded by $L^{-1}(eU_{2})$.

\section{\label{sec:EDL-varieties}Extended Deligne-Lusztig varieties}

As before, Let $F$ be an arbitrary local field with finite residue
field $\mathbb{F}_{q}$. Let $L_{0}$ be a finite totally ramified
Galois extension of $F$, and set $L=L_{0}^{\textrm{ur}}$. Then $L$
is a finite extension of $F^{\textrm{ur}}$ (cf.~\cite{ivan}, II~4),
and thus $L$ is a Henselian discrete valuation field with the same
residue field as $F^{\textrm{ur}}$, namely $\overline{\mathbb{F}}_{q}$.
We have the relation $\mathfrak{p}_{F}\mathcal{O}_{L}=\mathfrak{p}_{L}^{e}$,
where $e=[L_{0}:F]$ is the ramification index of $L_{0}/F$. 

Restriction of automorphisms gives a map\[
\alpha:\textrm{Gal}(L/F)\longrightarrow\textrm{Gal}(F^{\textrm{ur}}/F)\longiso\textrm{Gal}(\overline{\mathbb{F}}_{q}/\mathbb{F}_{q})\supset\mathbb{Z},\]
where the subgroup $\mathbb{Z}$ is generated by the Frobenius map
$x\mapsto x^{q}$. The corresponding Frobenius element in $\textrm{Gal}(F^{\textrm{ur}}/F)$
is denoted by $\phi_{F}$. Let $\Gamma=\Gamma(L/F)$ be the group
$\alpha^{-1}(\mathbb{Z})\subset\textrm{Gal}(L/F)$. This is a relative
variant of the Weil group and sits in the following commutative diagram.
 $$\SelectTips{cm}{10}\xymatrix@H=.6cm{1\ar[r] & \textrm{Gal}(L/F^{\textrm{ur}})\ar[r]\ar@{=}[d] & \Gamma(L/F)\ar[r]\ar@{_{(}->}[d] & \,\langle\phi_F\rangle\,\,\ar[r]\ar@{_{(}->}[d] & 1\\ 1\ar[r] & \textrm{Gal}(L/F^{\textrm{ur}})\ar[r] & \textrm{Gal}(L/F)\ar[r] & \textrm{Gal}(F^{\textrm{ur}}/F)\ar[r]\ar[d]^{\cong} & 1\\ &  &  & \textrm{Gal}(L/L_{0})\ar@{_{(}->}[lu]}$$

  We see that $\phi_{L_{0}}\in\text{Gal}(L/L_{0})$ defines an element
in $\Gamma$ which is not in $\Gal(L/F^{\mathrm{ur}})$. Hence $\Gamma$
is generated by $\textrm{Gal}(L/F^{\textrm{ur}})$ together with the
element $\varphi_{L_{0}}$. The group $\Gal(L_{0}/F)$ is naturally
isomorphic to $\Gal(L/F^{\mathrm{ur}})$, and we shall identify elements
in the former with their corresponding images in the latter.

From now on, let $\mathbf{G}$ be either $\GL_{n}$ or $\SL_{n}$,
viewed as group schemes over $\mathcal{O}_{F}$. Let $\mathbf{T}$
be the standard split maximal torus in $\mathbf{G}.$ Let $\mathbf{B}$
be the upper-triangular Borel subgroup scheme of $\mathbf{G}$, and
let $\mathbf{U}$ be the unipotent radical of $\mathbf{B}$.

Let $r\geq1$ be a natural number. Every automorphism $\sigma\in\textrm{Gal}(L/F)$
stabilises $\mathcal{O}_{L}$ and $\mathfrak{p}_{L}^{r}$, respectively
(cf.~\cite{ivan}, II~Lemma 4.1). Therefore, each $\sigma\in\textrm{Gal}(L/F)$
defines a morphism of $\mathcal{O}_{F}$-algebras $\sigma:\mathcal{O}_{L,r}\rightarrow\mathcal{O}_{L,r}$,
and hence a homomorphism of groups $\sigma:\mathbf{G}(\mathcal{O}_{L,r})\rightarrow\mathbf{G}(\mathcal{O}_{L,r})$.
Moreover, $\mathcal{O}_{L,r}$ has the structure of algebraic ring
(isomorphic to affine $r$-space over $\overline{\mathbb{F}}_{q}$),
and each $\sigma\in\Gamma$ such that $\sigma\in\alpha^{-1}(\mathbb{Z}_{\geq0})$
gives rise to an algebraic endomorphism of $\mathcal{O}_{L,r}$. Hence
each $\sigma\in\Gal(L/F^{\mathrm{ur}})$ and each non-negative power
of $\phi_{L_{0}}$ induces (via the canonical isomorphism $\mathbf{G}(\mathcal{O}_{L,r})\cong G_{L,r}$)
an endomorphism of the algebraic group $G_{L,r}$. For $\sigma\in\Gal(L/F^{\mathrm{ur}})$,
the resulting endomorphism of $G_{L,r}$ is also denoted by $\sigma$.
Furthermore, the Frobenius map $\phi_{L_{0}}\in\text{Gal}(L/L_{0})$
induces a Frobenius endomorphism of the algebraic group ${G}_{L,r}$,
which we denote by $\phi$. It is clear that $T_{L,r}$, $B_{L,r}$,
and $U_{L,r}$ are stable under $\phi$ and under each of the endomorphisms
induced by $\sigma\in\Gal(L/F^{\mathrm{ur}})$.

In Section~\ref{sec:The-unramified-approach} the finite group $G_{F,r}$
was identified with the fixed points of ${G}_{r}$ under a Frobenius
map. However, this is not the only way to realise $G_{F,r}$ as a
group of fixed points of a connected algebraic group. The following
lemma and its corollary make this more precise for tamely ramified
extensions. The following is an additive Hilbert 90 for powers of
the maximal ideal $\mathfrak{p}_{L}$.
\begin{lem}
\label{lem:Hilbert90}Suppose that $L_{0}/F$ is tamely ramified.
Then $\Gal(L_{0}/F)$ is cyclic. Let $\sigma$ be a generator of $\Gal(L_{0}/F)$,
$m\geq1$ be a natural number, and $y\in\mathfrak{p}_{L_{0}}^{m}$
be an element such that $\Tr_{L_{0}/F}(y)=0$. Then there exists an
element $x\in\mathfrak{p}_{L_{0}}^{m}$ such that $x-\sigma(x)=y$.\end{lem}
\selectlanguage{british}%
\begin{proof}
Since $L_{0}/F$ is totally and tamely ramified, the Galois group
$\Gal(L_{0}/F)$ is cyclic of order $e$ (cf.~\cite{ivan}, II 4.4).
\foreignlanguage{english}{Tamely ramified extensions are characterised
by the fact that $\Tr$ maps units to units. In particular $e=\Tr_{L_{0}/F}(1)$
is a unit in $\mathcal{O}_{L_{0}}$, and $\Tr_{L_{0}/F}(1/e)=1$.
Let\[
x=\sum_{n=1}^{e-1}\left(\sigma^{n}(1/e)\cdot\sum_{i=0}^{n-1}\sigma^{i}(y)\right).\]
Then $x\in\mathfrak{p}_{L}^{m}$, and it is easily verified that $x-\sigma(x)=y$.}\end{proof}
\selectlanguage{english}%
\begin{cor}
\label{cor:Gamma-fix}Suppose that $L_{0}/F$ is tamely ramified,
and let $r\geq1$ be a natural number. Then $\mathcal{O}_{L,r}^{\Gamma}=\mathcal{O}_{F,r'}$,
where $r'=[\frac{r-1}{e}]+1$.\end{cor}
\begin{proof}
Since $L_{0}/F$ is totally and tamely ramified, it is cyclic, and
we choose a generator $\sigma$ of $\Gal(L_{0}/F)$. Following our
convention, we also use $\sigma$ to denote the corresponding generator
of $\Gal(L/F^{\mathrm{ur}})$. Now $\Gamma$ is generated by $\phi_{L_{0}}$
and $\sigma$ and since $\mathcal{O}_{L,r}^{\phi_{L_{0}}}=\mathcal{O}_{L_{0},r}$,
it is enough to show that $\mathcal{O}_{L_{0},r}^{\sigma}=\mathcal{O}_{F,r'}$.
It is well-known that $(\mathfrak{p}_{L_{0}}^{r})^{\sigma}=\mathfrak{p}_{L_{0}}^{r}\cap\mathcal{O}_{F}=\mathfrak{p}_{F}^{r'}$,
where $r'=[\frac{r-1}{e}]+1$. The functor of $\sigma$-invariants
is left exact, so we have an injection $\mathcal{O}_{F,r'}=\mathcal{O}_{L_{0}}^{\sigma}/(\mathfrak{p}_{L_{0}}^{r})^{\sigma}\hookrightarrow\mathcal{O}_{L_{0},r}^{\sigma}$.
Lemma~\ref{lem:Hilbert90} shows that $H^{1}(L_{0}/F,\mathfrak{p}_{L_{0}}^{m})=0$,
and so this injection is surjective, and this yields the result. 
\end{proof}
Recall that a B\'e{}zout domain is an integral domain in which every
finitely generated ideal is principal. 
\selectlanguage{british}%
\begin{lem}
\label{lem:Triangulisation}Let $R$ be a \foreignlanguage{english}{B\'e{}zout
domain}, and let $x\in\GL_{n}(R)$ be an arbitrary element, where
$n\geq2$. Suppose that the characteristic polynomial of $x$ splits
into linear factors over $R$. Then there exists an element $\lambda\in\SL_{n}(R)$,
such that $\lambda^{-1}x\lambda\in\mathbf{B}(R)$.\end{lem}
\begin{proof}
Let $a_{1}\in R$ be an eigenvalue of $x$ with corresponding eigenvector
$v=\left(\begin{smallmatrix}v_{1}\\
\vdots\\
v_{n}\end{smallmatrix}\right)\in R^{n}$, so that $xv=a_{1}v$. If $g\in\GL_{n}(R)$, then $gv$ is obviously
an eigenvector of $g^{-1}xg$. We claim that we can choose $g$ such
that $gv$ has an entry equal to $1$. Without loss of generality,
we may assume that there exists two integers $1\leq m,m'\leq n$,
such that $\gcd(v_{m},v_{m'})=1$. Then, since $R$ is a \foreignlanguage{english}{B\'e{}zout
domain, there exist elements $\alpha,\beta\in R$ such that \[
\alpha v_{m}+\beta v_{m'}=1.\]
Let $g=(g_{ij})$ be the matrix such that $g_{mm}=\alpha$, $g_{mm'}=\beta$,
$g_{m'm}=-v_{m'}$, $g_{mm'}=v_{m}$, $g_{ii}=1$ for all $i\notin\{m,m'\}$,
and all other entries equal to $0$. We have $g\in\SL_{n}(R)$, and
the $m^{\mathrm{th}}$ entry of $gv$ equals $1$, which proves the
claim. This implies that there exists a matrix} $\lambda_{1}\in\SL_{n}(R)$
matrix whose first column is the vector $gv$. We then have \[
\lambda_{1}^{-1}g^{-1}xg\lambda_{1}=\begin{pmatrix}* & * & \cdots & *\\
0\\
\vdots &  & x_{1}\\
0\end{pmatrix},\]
where $x_{1}\in\GL_{n-1}(R)$. We can now repeat the process by choosing
an eigenvalue of $x_{1}$. Working inductively, we obtain an element
$\lambda\in\SL_{n}(R)$ such that $\lambda^{-1}x\lambda\in\mathbf{B}(R)$. 
\end{proof}
The above lemma shows in particular that for any $x\in\mathbf{G}(\mathcal{O}_{F^{\mathrm{ur}}})$,
there exists a finite field extension $L/F^{\mathrm{ur}}$, and an
element $\lambda\in\mathbf{G}(\mathcal{O}_{L})$ such that $\lambda^{-1}x\lambda\in\mathbf{B}(\mathcal{O}_{L})$.
Reducing modulo $\mathfrak{p}_{L}^{r}$ we see that for any $x\in G_{F,r'}$
with $r'$ such that $G_{F,r'}\subseteq G_{L,r}$, there exists a
$\lambda\in{G}_{L,r}$ such that $\lambda^{-1}x\lambda\in B_{L,r}$. 

Recall that an element $x\in G_{r}$ is called \emph{regular} if its
centraliser $C_{G_{r}}(x)$ has minimal dimension (cf.~\cite{Hill_regular}
or \cite{dignemichel}, 14). Note that this is a more general definition
than that given in \cite{borel}, 12.2 (which coincides with the notion
of regular semisimple).
\begin{defn}
An element in $\mathbf{G}(\mathcal{O}_{F^{\mathrm{ur}},r})$ is called
\emph{separable} if it has distinct eigenvalues. Similarly, an element
in $G_{r}$ is called \emph{separable} if its corresponding element
in $\mathbf{G}(\mathcal{O}_{F^{\mathrm{ur}},r})$ (via the canonical
isomorphism $G_{r}\cong\mathbf{G}(\mathcal{O}_{F^{\mathrm{ur}},r})$)
is separable. If $x\in G_{r}$ is a regular separable element, we
call its centraliser $C_{G_{r}}(x)$ a \emph{quasi-Cartan} subgroup
(of $G_{r}$). Similarly, we call the finite group $C_{G_{F,r}}(x)$
a \emph{quasi-Cartan} subgroup (of $G_{F,r}$).
\end{defn}
Note that if $r=1$, then an element is regular semisimple if and
only if it is separable. In general, regular semisimple elements in
$G_{r}$ are separable, but there also exist unipotent regular separable
elements.

From now on, let $x\in G_{r}$ be a regular separable element. Since
$x$ is regular we then have \[
C_{G_{r}}(x)=\mathcal{O}_{F^{\mathrm{ur}},r}[x]\cap G_{r}.\]
Let $L/F^{\mathrm{ur}}$ be a finite field extension and $r'\geq r$
a\foreignlanguage{english}{ natural number} such that ${G}_{r'}$
is a subgroup of ${G}_{L,r}$ and such that there exists an element
$\lambda\in{G}_{L,r}$ such that $\lambda^{-1}x\lambda\in{B}_{L,r}$
(which is possible thanks to Lemma~\ref{lem:Triangulisation}). From
now on, let \foreignlanguage{english}{$r'=[\frac{r-1}{e}]+1$. }Let
$\Sigma_{0}$ be a set of generators of the finite group $\Gal(L/F^{\mathrm{ur}})$,
and put $\Sigma:=\{\phi\}\cup\Sigma_{0}$. Notice that if $L_{0}/F$
is tamely ramified, then Lemma~\ref{lem:Hilbert90} and Corollary~\ref{cor:Gamma-fix}
show that we can take $\Sigma_{0}$ to be a one-element set, and that
${G}_{L,r}^{\Sigma}={G}_{L,r}^{\Gamma}=G_{r'}$.

A subgroup of ${G}_{L,r}$ conjugate to ${B}_{L,r}$ will be called
a \emph{strict Borel subgroup.} Strict Borel subgroups are solvable,
but are not in general Borel subgroups of the algebraic group ${G}_{L,r}$.
Since $x$ is regular, we see that the group $C_{{G}_{r}}(x)$ lies
in the strict Borel $\lambda{B}_{L,r}\lambda^{-1}$. 
\begin{lem}
\label{lem:Borel-self-norm}Assume that \foreignlanguage{english}{$\mathbf{G}$
is either $\GL_{n}$ or $\SL_{n}$}. Then strict Borel subgroups in
${G}_{L,r}$ are self-normalising, that is, if $g\in{G}_{L,r}$ and
$g{B}_{L,r}g^{-1}\subseteq{B}_{L,r}$, then $g\in{B}_{L,r}$.\end{lem}
\begin{proof}
It is sufficient to prove the assertion for the group ${B}_{L,r}$.
In \cite{lees}, Lemma 1.2, it is shown that $\mathbf{B}(R)$ is self-normalising
in $\GL_{n}(R)$, when $R$ is a finite local PIR. The same proof
goes through for rings of the form $\mathcal{O}_{L,r}$, so the assertion
holds for $\mathbf{G}=\GL_{n}$. Since for any ring $R$ we have $\GL_{n}(R)=Z(R)\SL_{n}(R)$,
where $Z(R)$ is the subgroup of scalar matrices, the corresponding
assertion for $\mathbf{G}=\SL_{n}$ follows. It remains to use the
isomorphisms $\mathbf{G}(\mathcal{O}_{L,r})\cong{G}_{L,r}$ and $\mathbf{B}(\mathcal{O}_{L,r})\cong{B}_{L,r}$.\end{proof}
\begin{lem}
\label{lem:L open-closed}Let $G$ be a connected algebraic group,
and $\phi:G\rightarrow G$ a Frobenius endomorphism, that is, $\phi$
is surjective and $G^{\phi}$ is finite. Then the corresponding Lang
map $L:G\rightarrow G$, $g\mapsto g^{-1}\phi(g)$ is an open and
closed morphism.\end{lem}
\begin{proof}
By the Lang-Steinberg theorem $L$ is surjective, so it is in particular
a dominant map of irreducible varieties. Let $W\subseteq G$ be a
closed irreducible subset. Since the fibres of $L$ are all of the
form $G^{\phi}x$, for $x\in G$, the map $L:L^{-1}(W)\rightarrow W$
is an orbit map. By \cite{borel}, II 6.4, $G^{\phi}$ then acts transitively
on the set of irreducible components of $L^{-1}(W)$, and hence they
all have the same dimension, equal to the dimension of $G^{\phi}\backslash L^{-1}(W)\cong W$.
By \cite{humphreys}, Theorem 4.5, the map $L$ is thus open.

Now let $X\subseteq G$ be a closed subset. The set $G^{\phi}X$ is
then a closed subset which is a union of fibres. Hence\[
L(G-G^{\phi}X)=L(G)-L(G^{\phi}X)=G-L(X),\]
and since $G-X$ is open, and $L$ is open, $L(X)$ is closed in $G$.
\end{proof}
Let $\mathcal{B}_{L,r}$ denote the set of strict Borel subgroups
of ${G}_{L,r}$. Since ${B}_{L,r}$ is self-normalising in ${G}_{L,r}$,
strict Borels are in one-to-one correspondence with points of the
variety ${X}_{L,r}:={G}_{L,r}/{B}_{L,r}$. Consider the product $\prod_{\sigma\in\{1\}\cup\Sigma}X_{L,r}$,
with ${G}_{L,r}$ acting diagonally. For $(B_{\sigma})_{\sigma\in\{1\}\cup\Sigma}\in\prod_{\sigma\in\{1\}\cup\Sigma}X_{L,r}$,
we thus have the corresponding ${G}_{L,r}$-orbit ${G}_{L,r}(B_{\sigma})_{\sigma\in\{1\}\cup\Sigma}$. 
\begin{defn}
We define the variety\begin{multline*}
X_{L,r}^{\Sigma}(\lambda)=\{B\in\mathcal{B}_{L,r}\mid{G}_{L,r}(\sigma(B))_{\sigma\in\{1\}\cup\Sigma}={G}_{L,r}(\sigma(\lambda{B}_{L,r}\lambda^{-1}))_{\sigma\in\{1\}\cup\Sigma}\}\\
=\{B\in\mathcal{B}_{L,r}\mid h(\sigma(B))_{\sigma\in\{1\}\cup\Sigma}=(\sigma(\lambda{B}_{L,r}\lambda^{-1}))_{\sigma\in\{1\}\cup\Sigma}\ \text{for some }h\in{G}_{L,r}\}.\end{multline*}

\end{defn}
Identifying $\mathcal{B}_{L,r}$ with ${X}_{L,r}$ we can rewrite
the variety as\begin{multline*}
X_{L,r}^{\Sigma}(\lambda)\\
=\{g\in{G}_{L,r}\mid\sigma(\lambda)^{-1}h\sigma(g)\in{B}_{L,r}\ \text{for all }\sigma\in\{1\}\cup\Sigma\text{ and some }h\in{G}_{L,r}\}/{B}_{L,r}\\
=\{g\in{G}_{L,r}\mid g^{-1}\sigma(g)\in b\lambda^{-1}\sigma(\lambda){B}_{L,r}\ \text{for all \ensuremath{\sigma\in\Sigma}}\text{ and some }b\in{B}_{L,r}\}/{B}_{L,r},\end{multline*}
and by making the substitution $g\mapsto gb^{-1}$, we can normalise
the defining relations so that \[
X_{L,r}^{\Sigma}(\lambda)=\{g\in{G}_{L,r}\mid g^{-1}\sigma(g)\in\lambda^{-1}\sigma(\lambda){B}_{L,r}\ \forall\,\sigma\in\Sigma\}/{B}_{L,r}(\lambda),\]
where \[
{B}_{L,r}(\lambda):=\bigcap_{\sigma\in\{1\}\cup\Sigma}\lambda^{-1}\sigma(\lambda){B}_{L,r}\sigma(\lambda)^{-1}\lambda.\]
From now on we will use this last model for $X_{L,r}^{\Sigma}(\lambda)$.
The finite group $G_{L,r}^{\Sigma}={G}_{L,r}^{\Gamma}$ acts on $X_{L,r}^{\Sigma}(\lambda)$
by left multiplication.

We would now like to define finite covers of the varieties $X_{L,r}^{\Sigma}(\lambda)$
in a way that naturally generalises the finite covers $\widetilde{X}_{r}(\hat{w})$,
defined in the unramified case where $L=F^{\mathrm{ur}}$, and $\hat{w}\in N_{{G}_{r}}({T}_{r})$.
In general, however, there does not seem to be any straightforward
way to define an analogous cover of the whole of $X_{L,r}^{\Sigma}(\lambda)$,
but only of a certain ${G}_{L,r}^{\Gamma}$-stable subvariety. For
ease of notation, write $\epsilon$ for $\lambda^{-1}\phi(\lambda)$.
Let \[
{A}:=\{\epsilon^{-1}b\epsilon\phi(b)^{-1}\mid b\in{B}_{L,r}(\lambda)\}.\]
Clearly, ${A}$ is the image of ${B}_{L,r}(\lambda)$ under the morphism
$G_{L,r}\rightarrow G_{L,r}$ given by the map $g\mapsto\epsilon^{-1}g\epsilon\phi(g)^{-1}$.
Thus $A$ is conjugate to the image of the map $g\mapsto g\epsilon\phi(g)^{-1}\epsilon^{-1}$,
which in turn is equal to the image of the map $g\mapsto g^{-1}\epsilon\phi(g)\epsilon^{-1}$.
This last map is the Lang map corresponding to the Frobenius endomorphism
$g\mapsto\epsilon\phi(g)\epsilon^{-1}$, so by Lemma~\ref{lem:L open-closed},
it sends $B_{L,r}(\lambda)$ to a closed set. Hence ${A}$ is a closed
subset of $G_{L,r}$. 

Define the following subvariety of $X_{L,r}^{\Sigma}(\lambda)$, given
by\[
X_{L,r}^{\Sigma}(\lambda,{A}):=\Big(\{g\in{G}_{L,r}\mid g^{-1}\phi(g)\in\epsilon{A}{U}_{L,r}\}\cap X_{L,r}^{\Sigma}(\lambda)\Big)/{B}_{L,r}(\lambda).\]
Note that ${B}_{L,r}(\lambda)$ acts on $\{g\in{G}_{L,r}\mid g^{-1}\phi(g)\in\epsilon{A}{U}_{L,r}\}$
by right multiplication, and that ${G}_{L,r}^{\Gamma}$ acts on $X_{L,r'}^{\Sigma}(\lambda,{A})$
by left multiplication. Since ${G}_{L,r}^{\Gamma}$ and ${B}_{L,r}(\lambda)$
act on $X_{L,r}^{\Sigma}(\lambda)$ and $X_{L,r}^{\Sigma}(\lambda,{A})$,
the complement $X_{L,r}^{\Sigma}(\lambda)\setminus X_{L,r}^{\Sigma}(\lambda,{A})$
is also stable under these actions. We can now normalise the defining
relations in $X_{L,r}^{\Sigma}(\lambda,{A})$ by using the action
of ${B}_{L,r}(\lambda)$, so that\[
X_{L,r}^{\Sigma}(\lambda,{A})=\Big(\{g\in{G}_{L,r}\mid g^{-1}\phi(g)\in\epsilon{U}_{L,r}\}\cap X_{L,r}^{\Sigma}(\lambda)\Big)/S(\lambda),\]
where\[
S(\lambda):=\{b\in{B}_{L,r}(\lambda)\mid\epsilon^{-1}b^{-1}\epsilon\phi(b)\in{U}_{L,r}\}.\]
Using the fact that ${B}_{L,r}(\lambda)\subseteq{B}_{L,r}$ normalises
${U}_{L,r}$, it is easy to see that $S(\lambda)$ is a subgroup of
${B}_{L,r}(\lambda)$. Moreover, $S(\lambda)$ contains ${U}_{L,r}\cap\epsilon{U}_{L,r}\epsilon^{-1}\cap{B}_{L,r}(\lambda)$
and acts on $\{g\in{G}_{L,r}\mid g^{-1}\phi(g)\in\epsilon{U}_{L,r}\}$
by right multiplication. Let $S(\lambda)^{0}$ denote the connected
component of $S(\lambda)$. We define the finite cover\[
\widetilde{X}_{L,r}^{\Sigma}(\lambda):=\Big(\{g\in{G}_{L,r}\mid g^{-1}\phi(g)\in\epsilon{U}_{L,r}\}\cap X_{L,r}^{\Sigma}(\lambda)\Big)/S(\lambda)^{0}\longrightarrow X_{L,r}^{\Sigma}(\lambda,{A}).\]
We see that the finite group $S(\lambda)/S(\lambda)^{0}$ acts on
$\widetilde{X}_{L,r}^{\Sigma}(\lambda)$. Together with the respective
${G}_{L,r}^{\Gamma}$-actions this clearly makes $\widetilde{X}_{L,r}^{\Sigma}(\lambda)\rightarrow X_{L,r}^{\Sigma}(\lambda,{A})$
a ${G}_{L,r}^{\Gamma}\times S(\lambda)/S(\lambda)^{0}$-equivariant
cover.
\selectlanguage{english}%
\begin{rem*}
We call the varieties $X_{L,r}^{\Sigma}(\lambda)$ and the covers
$\widetilde{X}_{L,r}^{\Sigma}(\lambda)$ \emph{extended} Deligne-Lusztig
varieties, for the following reasons. Firstly, the varieties typically
correspond to a (non-trivial) extension of the maximal unramified
extension. Secondly, the various groups involved are iterated extensions
of groups over the corresponding residue fields. Thirdly, there are
at least three other constructions which could be referred to as generalisations
of (certain) Deligne-Lusztig varieties, neither of which is in the
direction given here. One of these is the varieties of Deligne associated
to elements in certain braid monoids (cf.~\cite{deligne-braid});
another is the affine Deligne-Lusztig varieties of Kottwitz and Rapoport
(cf.~\cite{rapoport_affine_DL}), and the third is the varieties
of Digne and Michel \cite{non_connected_DL}, defined with respect
to not necessarily connected, reductive groups. 
\end{rem*}
\selectlanguage{british}%
We close this section by showing that extended Deligne-Lusztig varieties
are a natural generalisation of classical Deligne-Lusztig varieties
as well as of the varieties which appear in \cite{Lusztig-Fin-Rings}
and \cite{Alex_Unramified_reps} (in the case of general and special
linear groups over finite local PIRs with their standard Frobenius
maps $\phi$).

Let $ $$\mathbf{T}'$ be a maximal torus in \foreignlanguage{english}{$\mathbf{G}\times\mathcal{O}_{F^{\mathrm{ur}}}$}
such that the group ${T}_{r}'$ is $\phi$-stable. Then ${T}_{r}'=C_{{G}_{r}}(x)$,
for some regular semisimple element $x\in{G}_{r}^{\phi}$, and by
\cite{Alex_Unramified_reps},~2 we have ${T}_{r}'=\lambda{T}_{r}\lambda^{-1}$
for some $\lambda\in{G}_{r}$. Hence $\lambda$ is an element such
that $\lambda^{-1}x\lambda\in{T}_{r}\subseteq{B}_{r}$, and the condition
that ${T}_{r}'$ be $\phi$-stable implies that $\lambda^{-1}\phi(\lambda)\in N_{{G}_{r}}({T}_{r})$.
Let $\hat{w}:=\lambda^{-1}\phi(\lambda)$. Take $L_{0}=F$ (i.e.,
$L=F^{\mathrm{ur}}$), $r'=r$, so that $\Gamma=\langle\phi\rangle$,
and $\Sigma=\{\phi\}$. The resulting extended Deligne-Lusztig variety
attached to this data is\[
X_{F^{\mathrm{ur}},r}^{\{\phi\}}(\lambda)=\{g\in{G}_{r}\mid g^{-1}\phi(g)\in\hat{w}{B}_{r}\}/({B}_{r}\cap\hat{w}{B}_{r}\hat{w}^{-1}),\]
and since $\hat{w}$ normalises ${T}_{r}$ it follows that ${B}_{r}(\lambda)={T}_{r}({U}_{r}\cap\hat{w}{U}_{r}\hat{w}^{-1})$,
and the Lang-Steinberg theorem implies that ${A}\supseteq{T}_{r}$.
Hence $X_{F^{\mathrm{ur}},r}^{\{\phi\}}(\lambda,{A})=X_{F^{\mathrm{ur}},r}^{\{\phi\}}(\lambda)$.
Furthermore, we have\begin{multline*}
{S}(\lambda)=\{tu\in{T}_{r}({U}_{r}\cap\hat{w}{U}_{r}\hat{w}^{-1})\mid\hat{w}^{-1}u^{-1}t^{-1}\hat{w}\phi(tu)\in{U}_{r}\}\\
=\{tu\in{T}_{r}({U}_{r}\cap\hat{w}{U}_{r}\hat{w}^{-1})\mid\hat{w}t^{-1}\hat{w}\phi(t)\in{U}_{r}\}\\
=\{t\in{T}_{r}\mid\hat{w}t^{-1}\hat{w}\phi(t)=1\}({U}_{r}\cap\hat{w}{U}_{r}\hat{w}^{-1}),\end{multline*}
and so $S(\lambda)^{0}={U}_{r}\cap\hat{w}{U}_{r}\hat{w}^{-1}$ and
${S}(\lambda)/{S}(\lambda)^{0}\cong\{t\in{T}_{r}\mid\hat{w}t^{-1}\hat{w}\phi(t)=1\}$.
The corresponding cover is \[
\widetilde{X}_{F^{\mathrm{ur}},r}^{\{\phi\}}(\lambda)=\{g\in{G}_{r}\mid g^{-1}\phi(g)\in\hat{w}{U}_{r}\}/({U}_{r}\cap\hat{w}{U}_{r}\hat{w}^{-1}),\]
and hence $X_{F^{\mathrm{ur}},r}^{\{\phi\}}(\lambda)=X_{r}(\hat{w})$
and $\widetilde{X}_{F^{\mathrm{ur}},r}^{\{\phi\}}(\lambda)=\widetilde{X}_{r}(\hat{w})$
are the varieties we considered in Section~\ref{sec:The-unramified-approach}.
We thus see that the classical Deligne-Lusztig varieties as well as
the generalisations in \cite{Lusztig-Fin-Rings} and \cite{Alex_Unramified_reps}
(in the case of general or special linear groups over finite local
PIRs with their standard Frobenius maps $\phi$) appear as special
cases of the construction of extended Deligne-Lusztig varieties given
in this section.

\selectlanguage{english}%

\section{\label{sec:EDLforGL2-SL2}\foreignlanguage{british}{Extended Deligne-Lusztig
varieties for $\GL_{2}$ and $\SL_{2}$}}

\selectlanguage{british}%
Throughout this section $\mathbf{G}$ will denote either of the groups
$\GL_{2}$ or $\SL_{2}$, over $\mathcal{O}_{F}$. The subgroups $\mathbf{T}$,
$\mathbf{B}$, and $\mathbf{U}$ of $\mathbf{G}$ are the same as
in Section~\ref{sec:EDL-varieties}. As in the preceding section
we treat the two types of groups simultaneously in a uniform way.
Assume that $F$ is a local function field (i.e.,~$\chara F=p$).
Assume also that $F$ has residue characteristic different from $2$.
In this section we will study extended Deligne-Lusztig varieties for
groups of the form $G_{F,2}$. 

Let $\zeta$ denote an arbitrary fixed non-square unit in $\mathcal{O}_{F,2}$.
In $G_{F,2}$ the four distinct conjugacy classes of quasi-Cartans
are given by the following representatives:\begin{align*}
 & T_{F,2},\\
 & C_{G_{F,2}}\begin{pmatrix}0 & 1\\
\zeta & 0\end{pmatrix}=\left\{ \begin{pmatrix}a & b\\
\zeta b & a\end{pmatrix}\right\} \cap G_{F,2},\\
 & C_{G_{F,2}}\begin{pmatrix}0 & 1\\
\varpi & 0\end{pmatrix}=\left\{ \begin{pmatrix}a & b\\
\varpi b & a\end{pmatrix}\right\} \cap G_{F,2},\\
 & C_{G_{F,2}}\begin{pmatrix}0 & 1\\
\zeta\varpi & 0\end{pmatrix}=\left\{ \begin{pmatrix}a & b\\
\zeta\varpi b & a\end{pmatrix}\right\} \cap G_{F,2}.\end{align*}
The first two of these quasi-Cartans are \emph{unramified} in the
sense that each of them is the $\mathcal{O}_{F,2}$-points of some
maximal torus of the group scheme $\mathbf{G}$. They are also unramified
in the sense that they can be brought into triangular form over $\mathcal{O}_{F^{\mathrm{ur}},2}$,
that is, there exists a $\lambda\in{G}_{2}$ such that $\lambda^{-1}C_{G_{F,2}}\begin{pmatrix}0 & 1\\
\zeta & 0\end{pmatrix}\lambda\subseteq{B}_{2}$ (for $T_{F,2}$ this is a trivial fact). For the maximal torus $T_{F,2}$,
we can take $\lambda=1$, and this gives rise to the variety $X_{2}(1)$.
Each $\lambda$ that triangulises $C_{G_{F,2}}\begin{pmatrix}0 & 1\\
\zeta & 0\end{pmatrix}$ gives rise to the variety $X_{2}(\lambda)=X_{2}(\hat{w})$, where
$w$ is the non-trivial Weyl group element in $G_{1}$. Now the cover
$\widetilde{X}_{2}(\lambda)$ of $X_{2}(\lambda)$ depends on $\lambda$,
that is, on the choice of strict Borel subgroup containing the Cartan
subgroup in question. However, it is known that the possible finite
covers of $X_{2}(1)$ and $X_{2}(\hat{w})$ of the type we are considering
all give rise to equivalent representations $R_{\mathbf{T},\theta}$
in their cohomology (cf.~\cite{Alex_Unramified_reps}, Corollary~3.4).

We will refer to the last two of the above quasi-Cartans as \emph{ramified.}
We now attach extended Deligne-Lusztig varieties and corresponding
representations also to the ramified quasi-Cartans. Let $L_{0}=F(\sqrt{\varpi})$
be one of the two ramified quadratic extensions of $F$ (recall that
$p\neq2$, so we have only tame ramification). Then $L=L_{0}^{\mathrm{ur}}$
is independent of the choice of ramified quadratic extension of $F$.
The group $\Gamma$ is generated by the Frobenius $\phi_{L_{0}}$
together with an involution $\sigma\in\Gal(L/F^{\mathrm{ur}})$, so
we take $\Sigma=\{\phi,\sigma\}$. Let $r=3$, so that $\mathcal{O}_{L,3}^{\Gamma}=\mathcal{O}_{L,3}^{\Sigma}=\mathcal{O}_{F,2}$.
We then have $G_{L,3}^{\Gamma}=G_{F,2}$. Define the following elements
of $\mathbf{G}(\mathcal{O}_{L,3})$: \[
\lambda=\begin{pmatrix}1 & 0\\
\sqrt{\varpi} & 1\end{pmatrix},\qquad\mu=\begin{pmatrix}1 & 0\\
\sqrt{\zeta\varpi} & 1\end{pmatrix}.\]
Then we clearly have\[
\lambda^{-1}C_{G_{2}}\begin{pmatrix}0 & 1\\
\varpi & 0\end{pmatrix}\lambda\subseteq\mathbf{B}(\mathcal{O}_{L,3}),\qquad\mu^{-1}C_{G_{2}}\begin{pmatrix}0 & 1\\
\zeta\varpi & 0\end{pmatrix}\mu\subseteq\mathbf{B}(\mathcal{O}_{L,3}).\]
This defines the associated extended Deligne-Lusztig varieties\begin{align*}
 & X_{L,3}^{\Sigma}(\lambda)=\{g\in{G}_{L,3}\mid g^{-1}\phi(g)\in{B}_{L,3},\ g^{-1}\sigma(g)\in\lambda^{-1}\sigma(\lambda){B}_{L,3}\}/{B}_{L,3}(\lambda),\\
 & X_{L,3}^{\Sigma}(\mu)=\{g\in{G}_{L,3}\mid g^{-1}\phi(g)\in\mu^{-1}\phi(\mu){B}_{L,3},\ g^{-1}\sigma(g)\in\mu^{-1}\sigma(\mu){B}_{L,3}\}/{B}_{L,3}(\mu),\end{align*}
(note that $\phi(\lambda)=\lambda$, and that $\phi(\mu)=\sigma(\mu)=\mu^{-1}$).

The corresponding covers are given by\begin{align*}
 & \widetilde{X}_{L,3}^{\Sigma}(\lambda)=\{g\in{G}_{L,3}\mid g^{-1}\phi(g)\in{U}_{L,3},\ g^{-1}\sigma(g)\in\lambda^{-1}\sigma(\lambda){B}_{L,3}\}/{S}(\lambda)^{0},\\
 & \widetilde{X}_{L,3}^{\Sigma}(\mu)=\{g\in{G}_{L,3}\mid g^{-1}\phi(g)\in\mu^{-1}\phi(\mu){U}_{L,3},\ g^{-1}\sigma(g)\in\mu^{-1}\sigma(\mu){B}_{L,3}\}/{S}(\mu)^{0},\end{align*}
where\begin{align*}
{S}(\lambda) & =\{b\in{B}_{L,r}(\lambda)\mid b^{-1}\phi(b)\in{U}_{L,r}\},\\
{S}(\mu) & =\{b\in{B}_{L,r}(\lambda)\mid\phi(\mu)^{-1}\mu b^{-1}\mu^{-1}\phi(\mu)\phi(b)\in{U}_{L,r}\}.\end{align*}

\begin{thm}
\label{thm:EDL var GL2-SL2}Let $\mathbf{Z}$ be the centre of $\mathbf{G}$.
Then there exists a $G_{L,3}^{\Sigma}$-equivariant isomorphism\[
\widetilde{X}_{L,3}^{\Sigma}(\lambda)/(Z_{L,3}^{1})^{\phi}\cong G_{L,3}^{\Sigma}/(Z_{L,3}^{1})^{\Sigma}(U_{L,3}^{1})^{\Sigma}.\]
\end{thm}
\begin{proof}
We begin by determining ${S}(\lambda)$ explicitly. \foreignlanguage{english}{For
simplicity we shall write $e$ for $\lambda^{-1}\sigma(\lambda)$,
in what follows.} First consider ${B}_{L,r}(\lambda)={B}_{L,r}\cap e{B}_{L,r}e^{-1}$.
We write elements in $\mathcal{O}_{L,3}$ in the form $a_{0}+a_{1}\sqrt{\varpi}+a_{2}\varpi$,
where $a_{i}\in\overline{\mathbb{F}}_{q}$. We then have\begin{eqnarray*}
 &  & \phi(a_{0}+a_{1}\sqrt{\varpi}+a_{2}\varpi)=a_{0}^{q}+a_{1}^{q}\sqrt{\varpi}+a_{2}^{q}\varpi,\\
 &  & \sigma(a_{0}+a_{1}\sqrt{\varpi}+a_{2}\varpi)=a_{0}-a_{1}\sqrt{\varpi}+a_{2}\varpi.\end{eqnarray*}
Note in particular that $\phi$ and $\sigma$ commute. As usual, we
identify subgroups of $\mathbf{G}(\mathcal{O}_{L,3})$ with their
corresponding subgroups in ${G}_{L,3}$. Then\foreignlanguage{english}{\[
{B}_{L,r}(\lambda)=\left\{ \begin{pmatrix}a_{0}+a_{1}\sqrt{\varpi}+a_{2}\varpi & \frac{d_{1}-a_{1}}{2}+b_{1}\sqrt{\varpi}+b_{2}\varpi\\
0 & a_{0}+d_{1}\sqrt{\varpi}+d_{2}\varpi\end{pmatrix}\mid a_{i},b_{i}\in\overline{\mathbb{F}}_{q}\right\} \cap{G}_{L,r},\]
and so\[
{S}(\lambda)=\left\{ \begin{pmatrix}a_{0}+a_{1}\sqrt{\varpi}+a_{2}\varpi & \frac{d_{1}-a_{1}}{2}+b_{1}\sqrt{\varpi}+b_{2}\varpi\\
0 & a_{0}+d_{1}\sqrt{\varpi}+d_{2}\varpi\end{pmatrix}\mid a_{i}^{q}=a_{i},\ d_{i}^{q}=d_{i}\right\} \cap{G}_{L,r}.\]
Hence, the connected component of ${S}(\lambda)$ is \[
{S}(\lambda)^{0}={U}_{L,3}^{1},\]
and ${S}(\lambda)/{S}(\lambda)^{0}\cong Z_{L,1}^{\phi}(T_{L,3}^{1})^{\phi}=Z_{1}^{\phi}(T_{L,3}^{1})^{\phi}$.}

\selectlanguage{english}%
Let $Y:=\{g\in{G}_{L,3}\mid g^{-1}\phi(g)\in{U}_{L,3},\ g^{-1}\sigma(g)\in e{B}_{L,3}\}$,
so that $\widetilde{X}_{L,3}^{\Sigma}(\lambda)=Y/{U}_{L,3}^{1}$.
For $g\in Y$ we have $g^{-1}\phi(g)=u$, and $g^{-1}\sigma(g)=eb$,
for some $u\in{U}_{L,3}$, $b\in{B}_{L,3}$. The commutativity of
$\phi$ and $\sigma$ yields $\sigma(gu)=\phi(geb)$, and since $\phi(e)=e$
this implies \[
eb\sigma(u)=ue\phi(b).\]
Hence we obtain $e^{-1}ue\in{B}_{L,3}$, so that $u\in{U}_{L,3}\cap e{B}_{L,3}e^{-1}={U}_{L,3}^{1}$.
We thus have $Y=\{g\in{G}_{L,3}\mid g^{-1}\phi(g)\in{U}_{L,3}^{1},\ g^{-1}\sigma(g)\in e{B}_{L,3}\}$.
If we set \[
Y':=\{g\in{G}_{L,3}^{\phi}\mid g^{-1}\sigma(g)\in e{B}_{L,3}\}/(Z_{L,3}^{1})^{\phi}({U}_{L,3}^{1})^{\phi},\]
we then have a natural $G_{L,3}^{\Sigma}$-equivariant isomorphism\[
\widetilde{X}_{L,3}^{\Sigma}(\lambda)/(Z_{L,3}^{1})^{\phi}=Y/(Z_{L,3}^{1})^{\phi}{U}_{L,3}^{1}\longiso Y'.\]
Now the translation map $g\mapsto g\lambda^{-1}$ is an equivariant
isomorphism $Y'\iso Y'\lambda^{-1}$, and we have\[
Y'\lambda^{-1}=\{g\in{G}_{L,3}^{\phi}\mid g^{-1}\sigma(g)\in\sigma(\lambda){B}_{L,3}\sigma(\lambda)^{-1}\}/(Z_{L,3}^{1})^{\phi}\lambda({U}_{L,3}^{1})^{\phi}\lambda^{-1}.\]
If $g\in Y'\lambda^{-1}$, then $g^{-1}\sigma(g)\in\sigma(\lambda){B}_{L,3}\sigma(\lambda)^{-1}$,
and we then also have $g^{-1}\sigma(g)\in\lambda{B}_{L,3}\lambda^{-1}$,
since $\sigma$ has order $2$. Therefore $g^{-1}\sigma(g)\in\sigma(\lambda){B}_{L,3}\sigma(\lambda)^{-1}\cap\lambda{B}_{L,3}\lambda^{-1}$,
which is equivalent to \[
\lambda^{-1}g^{-1}\sigma(g)\lambda\in e{B}_{L,3}e^{-1}\cap{B}_{L,3}={B}_{L,3}(\lambda).\]
We thus have $g^{-1}\sigma(g)\in\lambda{B}_{L,3}(\lambda)\lambda^{-1}$.
Now, the image of the map $L_{\sigma}:{G}_{L,3}\rightarrow{G}_{L,3}$
given by $g\mapsto g^{-1}\sigma(g)$ clearly lies in ${G}_{L,3}^{1}$.
Thus\begin{multline*}
g^{-1}\sigma(g)\in\lambda{B}_{L,3}(\lambda)\lambda^{-1}\cap{G}_{L,3}^{1}\\
=\lambda\left\{ \begin{pmatrix}1+a_{1}\sqrt{\varpi}+a_{2}\varpi & b_{1}\sqrt{\varpi}+b_{2}\varpi\\
0 & 1+a_{1}\sqrt{\varpi}+d_{2}\varpi\end{pmatrix}\mid a_{i},b_{i}\in\overline{\mathbb{F}}_{q}\right\} \lambda^{-1}\cap{G}_{L,3}^{\phi},\end{multline*}
 and since $\lambda$ normalises the above set of matrices, we get\begin{multline*}
g^{-1}\sigma(g)\in\left\{ \begin{pmatrix}1+a_{1}\sqrt{\varpi}+a_{2}\varpi & b_{1}\sqrt{\varpi}+b_{2}\varpi\\
0 & 1+a_{1}\sqrt{\varpi}+d_{2}\varpi\end{pmatrix}\mid a_{i},b_{i}\in\overline{\mathbb{F}}_{q}\right\} \cap{G}_{L,3}^{\phi}\\
=(Z_{L,3}^{1})^{\phi}(T_{L,3}^{2})^{\phi}(U_{L,3}^{1})^{\phi}.\end{multline*}
Now we can obviously replace the relation $g^{-1}\sigma(g)\in(Z_{L,3}^{1})^{\phi}(T_{L,3}^{2})^{\phi}(U_{L,3}^{1})^{\phi}$
by $g^{-1}\sigma(g)\in(Z_{L,3}^{1})^{\phi}(T_{L,3}^{2})^{\phi}(U_{L,3}^{1})^{\phi}\cap L_{\sigma}({G}_{L,3}^{\phi})$,
without loss of generality. We thus have\begin{multline*}
Y'\lambda^{-1}\\
=\{g\in G_{L,3}^{\phi}\mid g^{-1}\sigma(g)\in(Z_{L,3}^{1})^{\phi}(T_{L,3}^{2})^{\phi}(U_{L,3}^{1})^{\phi}\cap L_{\sigma}({G}_{L,3}^{\phi})\}/(Z_{L,3}^{1})^{\phi}\lambda({U}_{L,3}^{1})^{\phi}\lambda^{-1}.\end{multline*}
One shows by direct computation that \[
L_{\sigma}((Z_{L,3}^{1})^{\phi}\lambda({U}_{L,3}^{1})^{\phi}\lambda^{-1})\supseteq(Z_{L,3}^{1})^{\phi}(T_{L,3}^{2})^{\phi}(U_{L,3}^{1})^{\phi}\cap L_{\sigma}({G}_{L,3}^{\phi}).\]
This implies that there is a natural equivariant isomorphism\[
Y'\lambda^{-1}\longiso{G}_{L,3}^{\Sigma}/((Z_{L,3}^{1})^{\phi}\lambda({U}_{L,3}^{1})^{\phi}\lambda^{-1})^{\Sigma}=G_{L,3}^{\Sigma}/(Z_{L,3}^{1})^{\Sigma}(U_{L,3}^{1})^{\Sigma}=G_{F,2}/Z_{F,2}^{1}U_{F,2}^{1}.\]
Since $\widetilde{X}_{L,3}^{\Sigma}(\lambda)/(Z_{L,3}^{1})^{\phi}\cong Y'\lambda^{-1}$,
the theorem is proved.
\end{proof}
\selectlanguage{english}%
The above theorem, together with \cite{dignemichel}, 10.10 (i) shows
that the variety $\widetilde{X}_{L,3}^{\Sigma}(\lambda)$ affords
the representation \[
\Ind_{Z_{F,2}^{1}U_{F,2}^{1}}^{G_{F,2}}\mathbf{1}\]
as a subrepresentation of its cohomology. In particular, for $\mathbf{G}=\SL_{2}$,
we have $Z_{F,2}^{1}=\{1\}$ (using $p\neq2$). Moreover, it is easy
to show that for $\mathbf{G}=\GL_{2}$, each nilpotent representation
of $\GL_{2}(\mathcal{O}_{F,2})$ is an irreducible constituent of
$\Ind_{B_{F,2}^{1}}^{G_{F,2}}\mathbf{1}$ (cf.~\cite{Hill_Jordan},
Lemma~2.12; note that we have defined nilpotent representations to
be primitive). Thus $\widetilde{X}_{L,3}^{\Sigma}(\lambda)$ affords
in particular all the nilpotent representations of $G_{F,2}$, both
for $\mathbf{G}=\SL_{2}$ and $\mathbf{G}=\GL_{2}$. Together with
the results of Lusztig \cite{Lusztig-Fin-Rings}, Section~3, this
proves that every irreducible representation of $\SL_{2}(\mathbb{F}_{q}[[\varpi]]/(\varpi^{2}))$,
with $p$ odd, appears in the cohomology of some extended Deligne-Lusztig
variety attached to a (possibly ramified) quasi-Cartan subgroup.

\section{Further directions}

In the proof of Theorem~\ref{thm:EDL var GL2-SL2}, the hypothesis
that $F$ be a function field was only used to calculate the explicit
form of the various groups involved, and the image of $L_{\sigma}$.
It is therefore likely that the argument can be extended to any non-archimedean
local field $F$ with $p\neq2$, using similar methods. Furthermore,
the question of whether the action of the finite group $S(\lambda)/S(\lambda)^{0}$
on $\widetilde{X}_{L,3}^{\Sigma}(\lambda)$ can be used to decompose
$\Ind_{Z_{F,2}^{1}U_{F,2}^{1}}^{G_{F,2}}\mathbf{1}$ into irreducible
components, remains open. However, the techniques used in the proof
of Theorem~\ref{thm:EDL var GL2-SL2} should prove useful for answering
this. Provided Lusztig's computations in \cite{Lusztig-Fin-Rings},
Section~3 could be carried out for $\GL_{2}$, it would follow from
the results of this paper that every irreducible representation of
$\GL_{2}(\mathbb{F}_{q}[[\varpi]]/(\varpi^{2}))$, with $p$ odd,
is realised by an extended Deligne-Lusztig variety.\bigskip

A natural problem is to generalise the construction of extended Deligne-Lusztig
varieties to reductive group schemes $\mathbf{G}$ over $\mathcal{O}_{F}$
other than $\GL_{n}$ or $\SL_{n}$. The ingredients required for
such a generalisation are as follows. First, one needs a generalisation
of Lemma~\ref{lem:Borel-self-norm} to any $\mathbf{G}$. This has
recently been given in \cite{Alex-RedGreenAlg}. Moreover, one would
need the result that any quasi-Cartan is contained in a strict Borel
subgroup of some $G_{L,r}$, which requires a version of Lemma~\ref{lem:Triangulisation}
for a Borel subgroup of $\mathbf{G}$.

It is also a natural question to ask whether our construction can
be extended to the wildly ramified case. When $L/F$ is tamely ramified,
we have shown that $G_{L,r}^{\Sigma}=G_{F,r'}$, but in the wildly
ramified case this may no longer hold. The difficulties in the wildly
ramified case are perhaps a reflection of the fact that the representation
theory of the $p$-adic group $\mathbf{G}(F)$ is radically different
in the wildly ramified case. In particular, one cannot expect in this
case that all the interesting representations are parametrised in
a straightforward way by data attached to maximal tori. Our present
construction can thus be seen as dealing efficiently only with the
cases where $L/F$ is tamely ramified. It should however be noted
that the only obstacle to defining extended Deligne-Lusztig varieties
in the wildly ramified case it due to the problem of descending from
$G_{L,r}$ to $G_{F,r'}$ by taking fixed-points. This is therefore
mainly a problem about Galois theoretic properties of finite ring
extensions. To go further in the wildly ramified case, it seems that
one has to consider either elements in $\Aut_{\mathcal{O}_{F,r'}}(\mathcal{O}_{L,r})$
other than those coming from elements in $\Gal(L/F)$, or a larger
field extension $E/L$, such that $E/F$ is tamely ramified.

A fundamental result of Deligne and Lusztig (cf.~\cite{delignelusztig},
Corollary~7.7) is that every irreducible representation of ${G}_{1}^{\phi}$
appears in the $l$-adic cohomology of some variety $\widetilde{X}_{1}(\hat{w})$.
An important question is whether something similar holds for the groups
$G_{r'}^{\phi}=G_{L,r}^{\Sigma}$, with respect to the extended Deligne-Lusztig
varieties $\widetilde{X}_{L,r}^{\Sigma}(\lambda)$. Some aspects of
the representation theory of the groups $\GL_{n}(\mathcal{O}_{F})$
are analogous to the representation theory of the $p$-adic group
$\GL_{n}(F)$. In particular, the construction of tamely ramified
supercuspidal representations via certain characters of maximal tori,
due to Howe \cite{Howe-Tame}, provides some of the motivation for
attaching extended Deligne-Lusztig varieties to quasi-Cartans. Given
this analogy, and the results obtained for nilpotent representations
in Section~\ref{sec:EDLforGL2-SL2}, we state the following open
problem: 
\begin{quote}
Suppose that $n$ is prime to $p$. Is it true that any irreducible
representation of $\GL_{n}(\mathcal{O}_{F,r'})$ which is a type for
a supercuspidal representation of $\GL_{n}(F)$, appears in the $l$-adic
cohomology of some extended Deligne-Lusztig variety $\widetilde{X}_{L,r}^{\Sigma}(\lambda)$?
\end{quote}
Here $r'=[\frac{r-1}{e}]+1$, with $e=e(L/F^{\mathrm{ur}})$, as before.
For the definition of types, see \cite{Henniart-appendix} and \cite{BK-types}.
In particular, any depth zero supercuspidal type on $\GL_{n}(\mathcal{O}_{F})$
factors through $\GL_{n}(k)$, corresponds to an unramified maximal
torus, and is realised in the cohomology of some variety $\widetilde{X}_{1}(\hat{w})$,
by the result of Deligne and Lusztig mentioned above. Moreover, the
results in Section~\ref{sec:EDLforGL2-SL2} show that every nilpotent
representation of $\GL_{2}(\mathcal{O}_{F,2})$, for $F$ a function
field, is realised by some $\widetilde{X}_{L,r}^{\Sigma}(\lambda)$.
Thus, the answer to the question is affirmative at least as far as
nilpotent types on $\GL_{2}(\mathcal{O}_{F,2})$ are concerned.\bigskip

It is interesting to ask about the possible connections between the
constructions in this paper, and the theory of character sheaves.
In \cite{Lusztig-char-sheaves-gen}, Lusztig discusses, among other
things, the possibility of defining character sheaves on $G_{r}$,
where $F$ is a function field, and $\mathbf{G}$ is a reductive group
scheme over $k_{F}$. The conjecture in \cite{Lusztig-char-sheaves-gen},
8 predicts that there is a theory of character sheaves on $G_{r}$
for generic principal series representations (i.e., those that correspond
to regular characters of a split unramified Cartan). However, Lusztig
remarks that one cannot expect to have a complete theory of character
sheaves on $G_{r}$, citing the irreducible representations of dimension
$q^{2}-1$ of $G_{F,2}$ (for $\mathbf{G}=\GL_{2}$, $F$ a function
field) as a reason for this. Note that these representations are nilpotent.
By the results in Section~\ref{sub:A-counterexample} for the closely
related case where $\mathbf{G}=\SL_{2}$, one may indeed expect that
the nilpotent representations cannot all be accounted for by character
sheaves on $G_{r}$. One of the principal aims of this paper has been
to demonstrate that the correct algebraic groups for constructing
nilpotent representations of $G_{r'}^{\phi}=G_{L,r}^{\Sigma}$ for
$\mathbf{G}=\GL_{2}$ or $\mathbf{G}=\SL_{2}$ in the tamely ramified
case, are not the {}``unramified'' groups $G_{r'}$, but groups
of the form $G_{L,r}$, where $L$ is a finite non-trivial extension
of $F^{\mathrm{ur}}$. One may therefore ask whether there exists
a theory of character sheaves on the groups $G_{L,r}$, pertaining
to (some of) the representations which do not correspond to character
sheaves on groups of the form $G_{r'}$.

\bibliographystyle{alex}
\bibliography{alex,Reps_over_fin_rings}
\selectlanguage{english}

\end{document}